\documentclass[11pt,draft]{amsart}

\usepackage[all]{xy}
\usepackage{amssymb}
\usepackage{mathrsfs}
\usepackage{cite}
\usepackage{amsfonts}
\newcommand{\rmv}[1]{}

\def\wh{\widehat}

\def\pv#1{\ensuremath{{\bf#1}}}

\def\inv{^{-1}}
\def\p{\varphi}

\def\J{\mathrel{{\mathscr J}}} 
\def\D{\mathrel{\mathscr D}} 

\def\e<{\leq _{E}}

\def\ov#1{\ensuremath{\overline {#1}}}

\def\malce{\mathbin{\hbox{$\bigcirc$\rlap{\kern-8.3pt\raise0,50pt\hbox{$\mathtt{m}$}}}}}

\def\1sk{^{(1)}}

\def\to{\rightarrow}

\def\Thmname{Theorem}
\def\Propname{Proposition}
\def\Lemmaname{Lemma}
\def\Definitionname{Definition}
\newtheorem{Thm}{\Thmname}[section]
\newtheorem{Prop}[Thm]{\Propname}
\newtheorem{Lemma}[Thm]{\Lemmaname}
{\theoremstyle{definition}
\newtheorem{Def}[Thm]{\Definitionname}}
{\theoremstyle{remark}
\newtheorem{Rmk}[Thm]{Remark}}
\newtheorem{Cor}[Thm]{Corollary}
{\theoremstyle{remark}
\newtheorem{Example}[Thm]{Example}}

{\theoremstyle{remark}
\newtheorem*{Claim}{Claim}}

\numberwithin{equation}{section}

\title{Strong Morita Equivalence of Inverse Semigroups}

\author{Benjamin Steinberg}
\address{School of Mathematics and Statistics\\
Carleton University \\
1125 Colonel By Drive\\
Ottawa, Ontario  K1S 5B6 \\
Canada}
\thanks{The author was supported in part by NSERC.  He also gratefully acknowledges the support of DFG}
\email{bsteinbg@math.carleton.ca}
\date{January 18, 2009}


\begin{document}
\begin{abstract}
We introduce strong Morita equivalence for inverse semigroups.  This notion encompasses Mark Lawson's concept of enlargement.  Strongly Morita equivalent inverse semigroups have Morita equivalent universal groupoids in the sense of Paterson and hence strongly Morita equivalent universal and reduced $C^*$-algebras.  As a consequence we obtain a new proof of a result of Khoshkam and Skandalis showing that the $C^*$-algebra of an $F$-inverse semigroup is strongly Morita equivalent to a cross product of a commutative $C^*$-algebra by a group.
\end{abstract}
\maketitle

\section{Introduction}
With the remarkable exception of finite group theory, in most areas of mathematics the equivalence relation of isomorphism is much too strong for the purposes of classification.  The aim of this article is to introduce a notion of strong Morita equivalence for inverse semigroups.

Our motivation comes from the theory of $C^*$-algebras.  As argued by Exel~\cite{Exel}, there is an emerging theory of ``combinatorial'' $C^*$-algebras that are constructed via the general pattern
\[\fbox{combinatorial object}\mapsto \fbox{inverse semigroup}\mapsto \fbox{groupoid}\mapsto \fbox{$C^*$-algebra}\] that is perhaps best exemplified by \[\fbox{alphabet}\mapsto \fbox{Cuntz semigroup}\mapsto \fbox{Cuntz groupoid}\mapsto \fbox{Cuntz algebra~\cite{cuntz}}\] or more generally by Cuntz-Krieger algebras~\cite{cuntzkrieger}, graph $C^*$-algebras and their variants~\cite{graphinverse,Paterson,higherrank,ultragraph,Exel}.  See~\cite{Exel} for further discussion.  (We remark that the Cuntz semigroup is known to semigroup theorists as the polycyclic inverse monoid~\cite{polycyclic}.)

Now, there is a well-known notion of Morita equivalence of groupoids, and moreover Morita equivalent groupoids have strongly Morita equivalent $C^*$-algebras~\cite{MRW,nonHausdorff,renaultdisintegration,muhlywilliams}. Often equivalence of groupoid algebras can be traced back to equivalence of their groupoids.  Similarly, cross product decompositions of groupoid algebras can frequently be traced back to semidirect product decompositions of the corresponding groupoids.  One would hope then that, for combinatorial $C^*$-algebras, these notions can all be transferred back to the inverse semigroup setting.  Of course, semidirect products of inverse semigroups are well established~\cite{Lawson}.  But until now a notion of Morita equivalence at the inverse semigroup level has been lacking.

The first hint that such a relationship should exist is a result of Khoshkam and Skandalis~\cite{Skandalis} that the $C^*$-algebra of an $F$-inverse semigroup is strongly Morita equivalent to a cross product of its maximal group image with a commutative $C^*$-algebra.  This is reminiscent of the result of Lawson~\cite{Lawson} that an $F$-inverse semigroup has an enlargement which is the semidirect product of its maximal group image with an idempotent commutative semigroup (i.e., a semilattice).  Moreover, the construction of the enlargement and the cross product decomposition both involve the notion of globalizing a partial group action~\cite{exelpartial,Lawsonkellendonk,abadie}. Let us recall here that Lawson~\cite{Lawsonenlarge,Lawson} defined an inverse semigroup $S$ to be an enlargement of an inverse subsemigroup $T$ if $STS=S$ and $TST=T$.  He argued that enlargement should be a special case of a more general theory of Morita equivalence that was yet to be developed.  Given Lawson's credo that enlargement is a special case of Morita equivalence, surely this must explain the result of Khoshkam and Skandalis.

This paper introduces the theory of Morita equivalence of inverse semigroups sought after by Lawson.  It encompasses the notion of enlargement as a special case.  Strong Morita equivalence of inverse semigroups implies equivalence of their associated groupoids (Theorem~\ref{Moritagroupoids}) and hence strong Morita equivalence of their $C^*$-algebras.  This yields a new proof of the result of Khoshkam and Skandalis on $F$-inverse semigroups.  We also show that strong Morita equivalence implies equivalence of classifying toposes~\cite{funktopos,topos} and hence isomorphism of cohomology groups~\cite{lausch,loganathan}.  A number of shared structural properties of strongly Morita equivalent inverse semigroups is established.

It should be mentioned that Talwar developed a theory of strong Morita equivalence for semigroups in general~\cite{Talwar1,Talwar2,Talwar3}.  Our notion of strong Morita equivalence implies Talwar's notion in the case of an inverse semigroup.  I do not believe the converse is true, although for inverse monoids both notions do coincide (and in fact they both coincide with enlargement in this case).

\section{Equivalence bimodules}
We begin by defining a \emph{Morita context} for inverse semigroups.  The notion is inspired by Rieffel's notion of an equivalence (or imprimitivity) bimodule~\cite{Rieffel,Rieffel2} and Talwar's strong Morita equivalence of semigroups~\cite{Talwar1,Talwar2}.  We first recall that an \emph{inverse semigroup} is a semigroup $S$ so that, for all $s\in S$, there exists a unique element $s^*\in S$ so that $ss^*s=s$ and $s^*ss^*=s^*$.  The reader is referred to Lawson's book~\cite{Lawson} for basic properties of inverse semigroups.  We shall frequently use that the idempotents of an inverse semigroup commute.  Also an inverse semigroup $S$ is partially ordered by $s\leq t$ if $s=et$ for some idempotent $e\in S$.  This ordering is compatible with the multiplication.  For idempotents, their product is their meet.

\begin{Def}[Morita context]\label{defMorita}
A \emph{Morita context} consists of a $5$-tuple $(S,T,X,\langle\ ,\ \rangle,[\ ,\ ])$ where $S$ and $T$ are inverse semigroups, $X$ is a set equipped with a left action by $S$ and a right action by $T$ that commute and \[\langle\ ,\ \rangle\colon X\times X\to S,\ [\ ,\ ]\colon X\times X\to T\] are surjective functions satisfying the following axioms (where $x,y,z\in X$, $s\in S$ and $t\in T$):
\begin{enumerate}
\item $\langle sx,y\rangle = s\langle x,y\rangle$;
\item $\langle y,x\rangle = \langle x,y\rangle^*$;
\item $\langle x,x\rangle x =x$;
\item $[x,yt] = [x,y]t$;
\item $[y,x] = [x,y]^*$;
\item $x[x,x] =x$;
\item $\langle x,y\rangle z = x[y,z]$.
\end{enumerate}
One calls $X$ (together with the two ``inner products'') an \emph{equivalence bimodule}.
\end{Def}

Note that the conjunction of (3) and (7) implies (6), and dually.
Two inverse semigroups $S$ and $T$ are said to be \emph{strongly Morita equivalent} if there exists an equivalence bimodule for them, i.e., there exists a Morita context $(S,T,X,\langle\ ,\ \rangle,[\ ,\ ])$.  We shall verify shortly that this is indeed an equivalence relation.

For example, recall that if $S$ is an inverse semigroup and $T$ is an inverse subsemigroup, then $S$ is termed by Lawson an \emph{enlargement} of $T$ if $STS=S$ and $TST=T$ (or equivalently $SE(T)S=S$ and $E(T)SE(T)=T$)~\cite{Lawson,Lawsonenlarge}.  For instance, if $e\in S$ is an idempotent such that $SeS=S$, then $S$ is an enlargement of $eSe$.  If $S$ is an inverse monoid with zero, $X$ is a set and $B_X(S)$ is the inverse semigroup consisting of all matrices of the form $sE_{ij}$, with $E_{ij}$ an $X\times X$-matrix unit, then $B_X(S)$ is an enlargement of $S$.   Lawson proved that if $S$ is an $F$-inverse semigroup, then $S$ has an enlargement of the form $X\rtimes G$ where $X$ is a semilattice containing $E(S)$ and $G$ is the maximal group image of $S$ (and conversely)~\cite{Lawson,Lawsonenlarge}.  It was suggested by Lawson~\cite{Lawson,Lawsonenlarge} that the notion of enlargement should be a special case of a more general theory of Morita equivalence.  The following proposition indicates that this article is the more general theory that Lawson was seeking.

\begin{Prop}\label{enlargement}
Let $S$ be an enlargement of $T$. Then $S$ and $T$ are strongly Morita equivalent.
\end{Prop}
\begin{proof}
Let $X=ST$ and define $\langle x,y\rangle = xy^*\in STS=S$, $[x,y] = x^*y\in TST=T$.  We verify that $X$ is an equivalence bimodule.  First we show that $\langle\ ,\ \rangle$ is surjective.  Let $s\in S$ and write $s=s_0ts_1$ with $s_0,s_1\in S$ and $t\in T$.  Then $s_0t,s_1^*t^*t\in ST=X$ and $\langle s_0t,s_1^*t^*t\rangle = s_0tt^*ts_1 = s$.
Clearly $\langle sx,y\rangle = sxy^*=s\langle x,y\rangle$.  Also $\langle y,x\rangle = yx^* = (xy^*)^*= \langle x,y\rangle^*$ and $\langle x,x\rangle x=xx^*x=x$.  The verifications for $[\ ,\ ]$ are similar.  Finally, we verify the compatibility between the inner products.  Indeed, $\langle x,y\rangle z = xy^*z = x[y,z]$.  This proves that $S$ and $T$ are strongly Morita equivalent.
\end{proof}

As a consequence of Proposition~\ref{enlargement}, strong Morita equivalence is a reflexive relation.  More generally, one easily shows that isomorphic inverse semigroups are strongly Morita equivalent.    Before verifying symmetry and transitivity it is convenient to establish some basic properties of Morita contexts.

\begin{Prop}\label{basicprops}
Let $(S,T,X,\langle\ ,\ \rangle,[\ ,\ ])$ be a Morita context.  Then, for $x,y\in X$ and $s\in S$, $t\in T$:
\begin{enumerate}
\item $\langle x,y\rangle \langle z,w\rangle = \langle x[y,z],w\rangle$;
\item $[x,y][z,w] = [x,\langle y,z\rangle w]$;
\item $\langle x,sy\rangle = \langle x,y\rangle s^*$;
\item $[xt,y] = t^*[x,y]$;
\item $\langle sx,sx\rangle = s\langle x,x\rangle s^*$;
\item $[xt,xt] = t^*[x,x]t$;
\item $\langle x,x\rangle\in E(S)$, $[x,x]\in E(T)$ for $x\in X$;
\item The maps $p(x)= \langle x,x\rangle$ and $q(x)= [x,x]$ give surjective maps $p\colon X\to E(S)$ and $q\colon X\to E(T)$, respectively;
\item $\langle xt,y\rangle = \langle x,yt^*\rangle$,  $\langle x,yt\rangle = \langle xt^*,y\rangle$;
\item $[sx,y] = [x,s^*y]$,  $[x,sy] = [s^*x,y]$.
\end{enumerate}
\end{Prop}
\begin{proof}
We just handle the case of $[\ ,\ ]$, as the other is dual.  First note that  \[[x,y][z,w] =[x,y[z,w]] = [x,\langle y,z\rangle w]\] establishing (2).
The computation \[[xt,y] = [y,xt]^*= ([y,x]t)^* = t^*[x,y]\] yields (4).  One immediately obtains (6) from (4).  Next we verify (7) by computing
$[x,x][x,x] = [x,x[x,x]]=[x,x]$ by Definition~\ref{defMorita}(6).  To prove (8), suppose that $e\in E(T)$ and $[x,y]=e$.  We claim $[ye,ye]=e$.  Indeed, \[[ye,ye] = e[y,y]e = e[y,y] = [x,y][y,y] =[x,y[y,y]] = [x,y]=e\] as required.

Finally, we prove (10). We first observe that
\begin{align*}
[sx,y][sx,y]^*[x,s^*y] &= [sx,y][y,sx][x,s^*y]\\  &
= [sx,y][y,\langle sx,x\rangle s^*y]\\ &=
 [sx,y][y,\langle sx,sx\rangle y] \\&=
[sx,y][y,sx][sx,y]\\&=
[sx,y]
\end{align*}
and so $[sx,y]\leq [x,s^*y]$.  Dually, $[s^*y,x]\leq [y,sx]$ and so, taking inverses, $[x,s^*y]\leq [sx,y]$. This establishes (10).
\end{proof}

From the proposition, it is clear that strong Morita equivalence is symmetric: just make $X$ a right $S$-set and a left-$T$ set using the inversion in $S$ and $T$;  one then switches the roles of $\langle\ ,\  \rangle$ and $[\ ,\ ]$.

As a consequence of Proposition~\ref{basicprops}, we obtain that each element $x\in X$, gives rise to homomorphisms $\epsilon_x\colon E(S)\to E(T)$ and $\eta_x\colon E(T)\to E(S)$.  These homomorphisms will play a role in constructing the Morita equivalence of the associated groupoids.

\begin{Prop}\label{createhom}
Fix $x\in X$ and let us define $\epsilon_x\colon E(S)\to E(T)$ and \mbox{$\eta_x\colon E(T)\to E(S)$} by $\epsilon_x(e) = [ex,ex]$ and $\eta_x(f) = \langle xf,xf\rangle$. Then $\epsilon_x$ and $\eta_x$ are homomorphisms.
\end{Prop}
\begin{proof}
We just handle the case of $\epsilon_x$, as that of $\eta_x$ is dual. Let \mbox{$e,f\in E(S)$}.  Using Proposition~\ref{basicprops}, we have
\begin{align*}
\epsilon_x(e)\epsilon_x(f) &= [ex,ex][fx,fx] = [ex,x][x,fx] = [ex,\langle x,x\rangle fx] \\ &= [ex,f\langle x,x\rangle x] = [ex,fx] = [ex,efx]=[efx,efx] = \epsilon_x(ef)
\end{align*}
as required.
\end{proof}

In order to verify the transitivity of strong Morita equivalence we must recall the definition of the tensor product a right $T$-set with a left $T$-set.  Details on the tensor product and homological algebra in this context can be found in~\cite{actsbook,Talwar2,Talwar3} for instance.

If $X$ is a right $T$-set, $Y$ is a left $T$-set and $Z$ is a set, then one says that a map $f\colon X\times Y\to Z$ is $T$-bilinear if $f(xt,y)=f(x,ty)$ for all $x\in X$, $y\in Y$ and $t\in T$.  The tensor product $X\otimes_T Y$ is the quotient of $X\times Y$ by the least equivalence relation $\sim$ such that $(xt,y)\sim (x,ty)$ for all $x\in X$, $y\in Y$ and $t\in T$.  The class of $(x,y)$ is denoted $x\otimes y$ and the map $X\times Y\to X\otimes_T Y$ given by $(x,y)\mapsto x\otimes y$ is the universal bilinear map.  It is easy to verify that if $S$ has a left action on $X$ and $U$ has a right action on $Y$ such that these actions commute with the actions of $T$, then $X\otimes_T Y$ admits well-defined commuting actions of $S$ and $U$ given by $s(x\otimes y) = sx\otimes y$ and $(x\otimes y)u= x\otimes yu$ where $s\in S$, $x\in X$, $y\in Y$ and $u\in U$.

We are now ready to prove that strong Morita equivalence is indeed an equivalence relation.

\begin{Prop}
Let $(S,T,X,\langle\ ,\ \rangle_S,[\ ,\ ]_T)$ and $(T,U,Y,\langle\ ,\ \rangle_T,[\ ,\ ]_U)$ be Morita contexts.  Then $(S,U,X\otimes_T Y,\langle\ ,\ \rangle,[\ ,\ ])$ is a Morita context where
\begin{align}\label{firstinnerproduct}
\langle x\otimes y,x'\otimes y'\rangle &= \langle x\langle y,y'\rangle_T,x'\rangle_S\\ \label{secondinnerproduct}
[x\otimes y,x'\otimes y'] &= [y,[x,x']_Ty']_U.
\end{align}
  Consequently, strong Morita equivalence is an equivalence relation.
\end{Prop}
\begin{proof}
First we need to verify that the ``inner products'' are well defined.  For each $(x,y)\in X\times Y$, we define a map $\theta_{x,y}\colon X\times Y\to S$ given by \[\theta_{x,y}(x',y')=\langle x\langle y,y'\rangle_T,x'\rangle_S.\]  Let us show that $\theta_{x,y}$ is $T$-bilinear.  Indeed, if $t\in T$, then $\theta_{x,y}(x't,y')=\langle x\langle y,y'\rangle_T,x't\rangle_S=\langle x\langle y,y'\rangle_T t^*,x'\rangle_S= \langle x\langle y,ty'\rangle_T,x'\rangle_S = \theta_{x,y}(x',ty')$ as required.  Thus we may now view $\theta_{x,y}$ as an element of $S^{X\otimes_T Y}$.  We have therefore defined a map $\theta\colon X\times Y\to S^{X\otimes_T Y}$ given by $\theta(x,y)=\theta_{x,y}$.  We claim that $\theta$ is $T$-bilinear. For if $t\in T$, then $\theta_{xt,y}(x',y') = \langle xt\langle y,y'\rangle_T,x'\rangle_S = \langle x\langle ty,y'\rangle_T,x'\rangle_S = \theta_{x,ty}(x',y')$.  Thus $\theta\colon X\otimes_T Y\to S^{X\otimes_T Y}$ and so we have a well-defined map $(X\otimes_T Y)\times (X\otimes_T Y)\to S$ given by \[(x\otimes y,x'\otimes y')\mapsto \theta_{x,y}(x',y') = \langle x\langle y,y'\rangle_T,x'\rangle_S,\] establishing that \eqref{firstinnerproduct} is well defined.  The case of \eqref{secondinnerproduct} is dual.

Next we check the surjectivity of $\langle\ ,\ \rangle$ and Axioms (1)--(3) of Definition~\ref{defMorita}; the corresponding facts for $[\ ,\ ]$ are of course dual.  Let $s\in S$.  We may write $s=\langle x,x'\rangle_S$ for some $x,x'\in X$.  Moveover, we may also write $[x,x]_T = \langle y,y'\rangle_T$ for some $y,y'\in Y$.  Then \[\langle x\otimes y,x'\otimes y'\rangle = \langle x\langle y,y'\rangle_T,x'\rangle_S = \langle x[x,x]_T,x'\rangle_S = \langle x,x'\rangle_S=s\] and so $\langle\ ,\ \rangle$ is surjective.    Next we turn to Axiom (1).  Here $\langle s(x\otimes y), x'\otimes y'\rangle = \langle sx\langle y,y'\rangle_T,x'\rangle_S=s\langle x\langle y,y'\rangle_T,x'\rangle_S =s\langle  x\otimes y, x'\otimes y'\rangle$.

For (2), we have that $\langle x'\otimes y',x\otimes y\rangle = \langle x'\langle y',y\rangle_T,x\rangle_S = \langle x',x\langle y,y'\rangle_T\rangle_S= \langle x\langle y,y'\rangle_T,x'\rangle_S^* = \langle x\otimes y,x'\otimes y'\rangle^*$. Axiom (3) follows since
\begin{align*}
\langle x\otimes y,x\otimes y\rangle (x\otimes y) &= \langle x\langle y,y\rangle_T,x\rangle_S(x\otimes y)\\ &= x\langle y,y\rangle_T[x,x]_T\otimes y\\  &= x[x,x]_T\otimes \langle y,y\rangle_T y\\ &= x\otimes y.
\end{align*}

It remains to verify the compatibility between the two inner products (Axiom (7)).  On the one hand
\begin{align*}
\langle x_1\otimes y_1,x_2\otimes y_2\rangle (x_3\otimes y_3) &= (\langle x_1\langle y_1,y_2\rangle_T,x_2\rangle_Sx_3)\otimes y_3 \\ &= x_1\langle y_1,y_2\rangle_T[x_2,x_3]_T\otimes y_3;
\end{align*}
 on the other hand
\begin{align*}
(x_1\otimes y_1)[x_2\otimes y_2,x_3\otimes y_3] &= x_1\otimes (y_1[y_2,[x_2,x_3]_Ty_3]_U) \\ &= x_1\otimes \langle y_1,y_2\rangle_T [x_2,x_3]_T y_3
\end{align*}
 yielding (7).  This completes the proof of the proposition.
\end{proof}

Our next proposition is an analogue of a result of Talwar~\cite{Talwar2,Talwar3}.

\begin{Prop}
Let $(S,T,X,\langle\ ,\ \rangle,[\ ,\ ])$ be a Morita context.  Make $X$ a left $T$-set by putting $tx=xt^*$.  Then $\langle\ ,\ \rangle\colon X\times X\to S$ is $T$-bilinear and induces a bijection $X\otimes_T X\to S$ given by $x\otimes y\mapsto \langle x,y\rangle$.  In fact, $S$ is isomorphic to the semigroup $X\otimes_T X$ with multiplication and inversion given by the formulas
\begin{align*}
(x\otimes y)(x'\otimes y') &= x[y,x']\otimes y'\\
(x\otimes y)^* &= y\otimes x,
\end{align*}
respectively.
A dual result holds for $T$.
\end{Prop}
\begin{proof}
First note that $\langle xt,y\rangle = \langle x,yt^*\rangle =\langle x,ty\rangle$, establishing the bilinearity.  Clearly the map $x\otimes y\mapsto \langle x,y\rangle$ is surjective by the definition of a Morita context.  Suppose that $\langle x,y\rangle =\langle x',y'\rangle$, whence $\langle y,x\rangle = \langle y',x'\rangle$.  Then
\begin{align*}
x\otimes y &= x[x,x]\otimes y = x\otimes y[x,x] = x\otimes \langle y,x\rangle x = x\otimes \langle y',x'\rangle x = x\otimes y'[x',x] \\ &= x[x,x']\otimes y' = \langle x,x\rangle x'\otimes y'[y',y'] = \langle x,x\rangle x'[y',y']\otimes y'\\ &= \langle x,x\rangle \langle x',y'\rangle y'\otimes y' = \langle x,x\rangle \langle x,y\rangle y'\otimes y' = \langle x,y\rangle y'\otimes y'\\ &= \langle x',y'\rangle y'\otimes y' = x'[y',y']\otimes y' = x'\otimes y',
\end{align*}
as required.  The multiplication and inverse rules now follow from the formulas  $\langle x,y\rangle\langle x',y'\rangle = \langle x[y,x'],y'\rangle$ and $\langle x,y\rangle^* = \langle y,x\rangle$.
\end{proof}

\section{\'Etale $S$-sets}
In order to prove that strong Morita equivalence of inverse semigroups yields an equivalence of their universal groupoids, we need to consider the category of \'etale $S$-sets, introduced in~\cite{topos} by the author and Funk in their work on the classifying topos of an inverse semigroup~\cite{funktopos}.

\begin{Def}[\'etale set]
If $S$ is an inverse semigroup, then an \emph{\'etale left $S$-set} is a pair $(X,p)$ where $X$ is a set equipped with a left action of $S$ on $X$ and $p\colon X\to E(S)$ is a map such that $p(x)x=x$ and $p(sx) = sp(x)s^*$, for all $x\in X$ and $s\in S$.  If $p$ is surjective, we say that $(X,p)$ has \emph{global support}.  A morphism $\p\colon (X,p)\to (Y,q)$ of \'etale left $S$-sets is an $S$-set morphism $\p\colon X\to Y$ such that the diagram
\[\xymatrix{X\ar[rr]^{\p}\ar[dr]_{p} &  &          Y\ar[dl]^{q}\\
  &       E(S)    &}\]
commutes.  \'Etale right $S$-sets are defined dually.  The category of \'etale right $S$-sets is denoted $\mathscr BS$ and is called the \emph{classifying topos} of $S$~\cite{topos}.  Of course, the categories of \'etale left $S$-sets and \'etale right $S$-sets are isomorphic using the inversion of $S$.
\end{Def}

For example, $S$ itself is an \'etale left set with global support with respect to the map $r\colon S\to E(S)$ given by $r(s) =ss^*$.  Indeed, $r(s)s = s$ and $r(st) = stt^*s^* = sr(t)s^*$.  In fact, $(S,r)$ is a torsion-free generator for the \'etendue $\mathscr BS$ (under the above isomorphism between left and right sets)~\cite{funktopos,topos}.  Also $(E(S),1_{E(S)})$ is an \'etale left $S$-set with respect to the action by conjugation and is the terminal object of this category.  One may easily verify that for a group $G$, $\mathscr BG$ is the category of right $G$-sets, whereas, for a semilattice $E$, the category $\mathscr BE$ is equivalent to the category of presheaves on $E$. It is routine to check that if $(X,p)$ is an \'etale left $S$-set, then $S$ acts by partial bijections on $X$ by restricting the domain of $s\in S$ to $s^*sX$ and the codomain to $ss^*X$.  Applying this to $(S,r)$ one obtains the famous Preston-Wagner representation~\cite{Lawson}, whereas applying this to $(E(S),1_{E(S)})$ yields the Munn representation~\cite{Lawson}.

We shall say that inverse semigroups $S$ and $T$ are \emph{Morita equivalent} if $\mathscr BS$ and $\mathscr BT$ are equivalent categories.  It was shown by Funk that if $S$ is an enlargement of $T$, then the inclusion of $T$ into $S$ induces an equivalence of $\mathscr BT$ and $\mathscr BS$ and hence $S$ and $T$ are Morita equivalent~\cite{funktopos}.  We shall see later that strong Morita equivalence implies Morita equivalence.  Unfortunately, we do not know whether the converse holds.

A key property of \'etale $S$-sets is that they come equipped with a natural partial order preserved by the action of $S$.

\begin{Prop}\label{etalesets}
Let $X$ be an \'etale left $S$-set.  Then, for $x,y\in X$, the following are equivalent:
\begin{enumerate}
\item $x= p(x)y$;
\item $x=ey$ some $e\in E(S)$.
\end{enumerate}
Moreover, if we define $x\leq y$ if $x=p(x)y$, then $\leq$ is a partial order compatible with the left action of $S$ and $p$ is order preserving.  A dual result holds for \'etale right sets.
\end{Prop}
\begin{proof}
Obviously (1) implies (2).  For (2) implies (1), observe that $p(x)y = p(ey)y = ep(y)ey = ep(y)y=ey=x$.  Clearly $x=p(x)x$, yielding reflexivity.  If $x\leq y$ and $y\leq z$, then $x=p(x)y=p(x)p(y)z$ and so by (2) $x\leq z$.  Next, if $x\leq y$ and $y\leq x$, then $x = p(x)y = p(x)p(y)x = p(y)p(x)x = p(y)x=y$.  To see compatibility of the action, if $x\leq y$ and $s\in S$, then $p(sx)sy = sp(x)s^*sy = sp(x)y = sx$ and so $sx\leq sy$.  Finally, we check that $p$ is order preserving.  If $x=p(x)y$, then $p(x) = p(p(x)y) = p(x)p(y)p(x)\leq p(y)$.
\end{proof}

\begin{Rmk}
The reason for the terminology \'etale is that if we make posets into topological spaces by taking the downsets as the open sets, then, for $(X,p)$ an \'etale left $S$-set, the map $p\colon X\to E(S)$ is \'etale in the usual topological sense.
\end{Rmk}

If $S$ is an inverse semigroup, then there is a group $G(S)$ and a surjective homomorphism $\sigma\colon S\to G(S)$ so that every homomorphism from $S$ to a group factors through $G(S)$.  One calls $G(S)$ the \emph{maximal group image} of $S$.  It is defined as the quotient of $S$ by the germ equivalence relation, which identifies two elements of $S$ if they have a common lower bound in the natural partial order on $S$~\cite{Lawson}.  More generally,
every \'etale $S$-set $X$ has a set of germs $\varinjlim X$, which turns out to be a $G(S)$-set.  We shall need the following lemma to construct the germ equivalence relation.

\begin{Lemma}\label{meetsinetale}
Let $(X,p)$ be an \'etale left $S$-set and suppose that $y,z\in X$ have a common upper bound $x\in X$.  Then $p(y)z=p(z)y$ is the meet of $y,z$ in $X$.
\end{Lemma}
\begin{proof}
First note that $p(y)z\leq z$ and $p(z)y\leq y$, so our first goal is to prove the equality.  Indeed, $p(y)z = p(y)p(z)x =p(z)p(y)x=p(z)y$ by the definition of the ordering. Set $u=p(y)z=p(z)y$ and suppose that $w\leq y,z$.  Then we have
\begin{align*}
p(w)u &= p(w)p(y)z = p(y)p(w)z=p(y)w\\&=p(y)p(w)y=p(w)p(y)y=p(w)y=w
\end{align*}
and hence $w\leq u$.
\end{proof}

It is now immediate that the relation of having a common lower bound is an equivalence relation.

\begin{Prop}\label{germsofetalesets}
Let $(X,p)$ be an \'etale left $S$-set.  Define $x\sim y$ if there exists $z\in X$ with $z\leq x,y$.  Then $\sim$ is an equivalence relation and $\varinjlim X=X/{\sim}$ has the natural structure of a left $G(S)$-set defined by $\sigma(s)[x]= [sx]$.
\end{Prop}
\begin{proof}
Clearly $\sim$ is reflexive and symmetric.  Assume that $x\sim y$ and $y\sim z$.  Let $u\leq x,y$ and $v\leq y,z$.  Because $y$ is a common upper bound on $u,v$, they have a meet $w$ by Lemma~\ref{meetsinetale}.   Clearly $w\leq x,z$.  Thus $\sim$ is an equivalence relation.

Suppose that $\sigma(s)=\sigma(t)$ and $[x]=[y]$.  Then we can find $r\leq s,t$ and $w\leq x,y$.  Clearly $rw\leq sx,ty$.  Thus $[sx]=[ty]$ and the action of $G(S)$ is well defined.  It is then trivial to verify we have defined an action.
\end{proof}

In fact, what we have really done is defined a functor $\varinjlim\colon \mathscr BS\to \mathscr BG(S)$.   It is shown in~\cite{topos} that this functor restricts to an equivalence between the locally constant objects of $\mathscr BS$ and the category $\mathscr BG(S)$ and hence $G(S)$ is the fundamental group of the topos $\mathscr BS$, and in particular only depends on the Morita equivalence class of $S$.

The reason for this digression into \'etale sets is that a Morita context gives rise to compatible \'etale set structures on the equivalence bimodule.  The partial ordering on the equivalence bimodule plays an important role in constructing the Morita equivalence of the universal groupoids.

\begin{Prop}\label{etalestruct}
Let $(S,T,X,\langle\ ,\ \rangle,[\ ,\ ])$ be a Morita context.  Define $p\colon X\to E(S)$ by $p(x) = \langle x,x\rangle$ and $q\colon X\to E(T)$ by $q(x) = [x,x]$.  Then $(X,p)$ is an \'etale left $S$-set with global support, $(X,q)$ is an \'etale right $S$-set with global support.  Moreover, the order on $X$ induced by the structures $(X,p)$ and $(X,q)$ coincide.    In particular, the set of germs of $X$ is independent of which \'etale structure we put on $X$.  In addition, $x\leq x',y\leq y'$ implies $\langle x,y\rangle \leq \langle x',y'\rangle$ and $[x,y]\leq [x',y']$.
\end{Prop}
\begin{proof}
Proposition~\ref{basicprops} implies everything except the coincidence of the two partial orders and the monotonicity of the bilinear forms.   Let us write $x\leq _S y$ if $x=ey$ with $e\in E(S)$ and $x\leq_T y$ if $x=yf$ with $f\in E(T)$.  Suppose $x\leq_S y$; say $x=ey$ with $e\in E(S)$.  Then $y[x,x] = y[ey,ey] = \langle y,ey\rangle ey = \langle y,y\rangle ey =e\langle y,y\rangle y =x$. Thus $x\leq_T y$.  Dually, $x\leq_T y$ implies $x\leq_S y$ and so the two orders coincide.  Hence the two germ constructions coincide.  We verify the final statement only for $[\ ,\ ]$, as the other case is identical.  Indeed, $[x,y] = [x'q(x),y'q(y)] = q(x)[x',y']q(y)\leq [x',y']$ as required.
\end{proof}

\section{Morita equivalence of groupoids and $C^*$-algebras}
In this section, all inverse semigroups are assumed to be countable and discrete as in~\cite{Paterson}.
Paterson associated to each inverse semigroup an \'etale groupoid $\mathscr G(S)$, called its \emph{universal groupoid}, with the property that both the respective universal and reduced $C^*$-algebras of $S$ and $\mathscr G(S)$ coincide~\cite{Paterson}.   There is a well-known notion of Morita equivalence for locally compact groupoids~\cite{MRW,nonHausdorff}.  A celebrated theorem of Muhly, Renault and Williams~\cite{MRW} shows that Morita equivalent Hausdorff groupoids have strongly Morita equivalent universal and reduced $C^*$-algebras.  In general, $\mathscr{G}(S)$ is not Hausdorff.  However, it is a special case of a result of
Renault~\cite[Corollaire 5.4]{renaultdisintegration} that also in the non-Hausdorff setting the strong Morita equivalence of universal $C^*$-algebras remains valid for Morita equivalent groupoids\footnote{The author thanks Jean Renault for pointing him toward this result.}; see also~\cite{muhlywilliams}.   The corresponding result for reduced $C^*$-algebras in the non-Hausdorff context was established by Tu~\cite{nonHausdorff}.

Let us recall the necessary definitions.
By a \emph{locally compact space}, we mean a topological space $Z$ with a basis of locally compact Hausdorff spaces.  We do not assume that $Z$ itself is Hausdorff.
Since we shall only be interested in \'etale groupoids, we stick with this setting for the definition of Morita equivalence.  An \emph{\'etale groupoid} (also known as an $r$-discrete groupoid~\cite{Renault,Paterson,Exel}) is a topological groupoid $\mathscr G$ such that the unit space $\mathscr G^0$ is locally compact Hausdorff and the domain map $d$ (or equivalently the range map $r$) is an \'etale map, i.e., a local homeomorphism.   In general, $\mathscr G$ need not be Hausdorff.  As usual  we do not distinguish between the object set of a groupoid and its unit space.

If $\mathscr G$ is a topological groupoid, $Z$ is a locally compact space and $f\colon Z\to \mathscr G^0$ is a continuous map, then we can form a topological groupoid $\mathscr G[Z]$ with unit space $Z$ and whose arrows are of the form $(z',g,z)$ such that $z',z\in Z$ and $d(g)=f(z), r(g)=f(z')$ equipped with the subspace topology.   The product is given by $(z'',g',z')(z',g,z) = (z'',g'g,z)$ and the inverse by $(z',g,z)\inv = (z,g\inv,z')$.   The topology is induced by the product topology on $Z\times \mathscr G\times Z$.  One should think of $\mathscr G[Z]$ as an amplification of $\mathscr G$ by $Z$.  Notice that the unit space of $\mathscr G[Z]$ need not be Hausdorff as we have made no separability assumptions on $Z$.

\begin{Def}[Morita equivalence]
Let $\mathscr G_1$ and $\mathscr G_2$ be \'etale groupoids.  Then $\mathscr G_1$ and $\mathscr G_2$ are said to be \emph{Morita equivalent} if there is a locally compact space $Z$ and continuous open surjections $f_i\colon Z\to \mathscr G_i^0$, $i=1,2$, such that the groupoids $\mathscr G_1[Z]$ and $\mathscr G_2[Z]$ are topologically isomorphic.
\end{Def}

Tu observes in~\cite[Proposition 2.29]{nonHausdorff} that one may always assume that $Z$ is in fact Hausdorff.    There is an equivalent formulation in terms of ``bimodules'' that is probably more aesthetically appealing, but we shall not use it; see~\cite[Proposition 2.29]{nonHausdorff}.

Next we recall Paterson's construction of the universal groupoid of an inverse semigroup~\cite{Paterson,Exel}. See also~\cite{resendeetale} for more on the relationship between \'etale groupoids and inverse semigroups.  Let $S$ be an inverse semigroup with idempotent set $E$.  Let $\wh E$ be the space of non-zero homomorphisms \mbox{$\p\colon E\to \{0,1\}$} with the topology of pointwise convergence.  Define, for $e\in E$, the compact open set $D(e) =\{\p\in \wh E\mid \p(e)=1\}$.  The inverse semigroup $S$ has a natural action on $\wh E$ as a pseudogroup of transformations.  If $s\in S$, then the domain of the action of $S$ is $D(s^*s)$ and the range is $D(ss^*)$.  The action is given by $s\p(e) = \p(s^*es)$ for $\p \in D(s^*s)$.

The universal groupoid of $S$ is the groupoid of germs for this action; the unit space of $\mathscr G(S)^0$ is $\wh E$ and arrows are equivalence classes $\ov {(s,\p)}$ where $\p\in D(s^*s)$ and $(s,\p)$ is equivalent to $(s',\psi)$ if and only if $\p=\psi$ and there exists $u\in S$ so that $u\leq s,s'$ and $\p\in D(u^*u)$.   The topology has a basis consisting of all sets of the form $(s,U)$ where $U\subseteq D(s^*s)$ is an open subset of $\wh E$ and $(s,U) = \{\ov{(s,\p)}\mid \p\in U\}$.  By definition, $d(\ov{(s,\p)}) = \p$ and $r(\ov{(s,\p)} = s\p$.   Multiplication of composable arrows is given by $\ov{(s,\p)}\cdot\ov{(s',\psi)} = \ov {(ss',\psi)}$ and inversion by $\ov{(s,\p)}\inv = \ov{(s^*,s\p)}$.  One can verify that $\mathscr G(S)$ is an \'etale groupoid~\cite{Exel,Paterson}.  It is known that if $S$ is $E$-unitary (or more generally $0$-$E$-unitary), then $\mathscr G(S)$ is Hausdorff~\cite{Paterson,Exel}; recall that an inverse semigroup is called $E$-unitary if any element above an idempotent is an idempotent~\cite{Lawson}.

The reader is referred to Proposition~\ref{createhom} for the definitions of $\eta_x$ and $\epsilon_x$ in the following proposition.

\begin{Prop}\label{homeomorphisms}
Let $(S,T,X,\langle\ ,\ \rangle,[\ ,\ ])$ be a Morita context and fix $x\in X$.  Then we can define inverse homeomorphisms $\alpha_x\colon D([x,x])\to D(\langle x,x\rangle)$ and $\beta_x\colon D(\langle x,x\rangle)\to D([x,x])$ by $\alpha_x(\p) = \p\circ \epsilon_x$ and $\beta_x(\psi) = \psi\circ \eta_x$.
\end{Prop}
\begin{proof}
By Proposition~\ref{createhom}, $\epsilon_x$ and $\eta_x$ are homomorphisms, and so they induce continuous maps $\alpha_x\colon \wh{E(T)}\to \wh{E(S)}$ and \mbox{$\beta_x\colon \wh{E(S)}\to \wh {E(T)}$} via $\alpha_x(\p) = \p\circ \epsilon_x$ and $\beta_x(\psi) = \psi\circ \eta_x$ by functoriality.   Thus we just need to check that $\alpha_x$ and $\beta_x$ are inverses when restricted to the above open subsets.  Suppose $\p \in D([x,x])$ and $e\in E(T)$.  Then
\[(\beta_x(\alpha_x(\p)))(e) = \alpha_x(\p)(\langle xe,xe\rangle) = \p([\langle xe,xe\rangle x,\langle xe,xe\rangle x]).\]  But $\langle xe,xe\rangle x = xe[xe,x] = xe[x,x] = x[x,x]e = xe$ and so \[(\beta_x(\alpha_x(\p)))(e) = \p([xe,xe]) = \p(e[x,x]e) = \p([x,x])\p(e) = \p(e)\] since $\p\in D([x,x])$.  The proof that $\alpha_x\beta_x$ is the identity is symmetric.
\end{proof}

Frequently, we write $x\p$ for $\alpha_x(\p)$ and $x^*\psi$ for $\beta_x(\psi)$.  The next lemma, establishing the compatibility of the actions of $S$, $T$ and $X$, will often be used without explicit reference in the sequel.

\begin{Lemma}\label{explicitref}
Let $x\in X$ and $s\in S$. Then $\p\in D([x,x])\cap x^*D(s^*s)$ if and only if $\p\in D([x,x])\cap D([sx,sx])$.  In this case, $(sx)\p = s(x\p)$.  A dual result holds for elements of $T$.
\end{Lemma}
\begin{proof}
Suppose $\p\in D([x,x])$.  Then $x\p(s^*s) =\p([s^*sx,s^*sx]) = \p([sx,sx])$ (using Proposition~\ref{basicprops}) and so $\p\in D([sx,sx])$ if and only if $\p\in x^*D(s^*s)$ by Proposition~\ref{homeomorphisms}.   To see that $(sx)\p = s(x\p)$, let $e\in E(S)$.  Then we have $((sx)\p)(e) = \p([esx,esx]) = \p([ss^*esx,esx]) =\p([s^*esx,s^*esx]) =x\p(s^*es) = (s(x\p))(e)$ as required.
\end{proof}

Next we want to show that the ordering on $X$ coming from Proposition~\ref{etalestruct} corresponds with restriction.

\begin{Prop}\label{restriction}
Let $(S,T,X,\langle\ ,\ \rangle,[\ ,\ ])$ be a Morita context.   Suppose that $x,y\in X$ with $x\leq y$.  Then $D([x,x])\subseteq D([y,y])$ and $\alpha_x$ is the restriction of $\alpha_y$.  The analogous result holds for $\beta_x,\beta_y$.
\end{Prop}
\begin{proof}
As $x\leq y$, we have $x=y[x,x]$ and so if $\p\in D([x,x])$, then \[1=\p([x,x]) = \p([y[x,x],y[x,x]]) = \p([x,x][y,y][x,x]) = \p([y,y]).\] Thus $\p\in D([y,y])$.  Assume now that $e\in E(S)$. Then we have $[ex,ex] = [ey[x,x],ey[x,x]] = [x,x][ey,ey][x,x]$.  Since $\p([x,x])=1$, it follows that $x\p(e)=\p([ex,ex])=\p([ey,ey]) = y\p(e)$.  This completes the proof.
\end{proof}

Let $(S,T,X,\langle\ ,\ \rangle,[\ ,\ ])$ be a Morita context.   We now define our space $Z$ as the space of germs of the $\alpha_x$ with $x\in X$.  Formally, $Z$ consists of all equivalence classes $\ov{(x,\p)}$ where $x\in X$ and $\p\in D([x,x])$.  Here $(x,\p)\sim(x',\psi)$ if $\p=\psi$ and there exists $y\leq x,x'$ with $\p\in D([y,y])$.  We take as a basis for a topology on $Z$ all sets of the form $(x,U)$ with $x\in X$ and $U\subseteq D([x,x])$ an open subset, where $(x,U) = \{\ov{(x,\p)}\mid \p\in U\}$.  Notice that the definition of $Z$ seems asymmetric in that $T$ appears in the definition but not $S$.  But this is an illusion: we could alternatively define $Z'$ to consist of equivalence classes pairs $\ov{(x,\p)}$ with $\p \in D(\langle x,x\rangle)$ and with a basis for the topology consisting of open sets $(x,U)$ where $U\subseteq D(\langle x,x\rangle)$ is open in $\wh{E(S)}$ and $(x,U) = \{\ov{(x,\p)}\mid \p\in U\}$.  Then $Z'$ is homeomorphic to $Z$ via the map $\ov{(x,\p)}\mapsto \ov{(x,x^*\p)}$ thanks to Proposition~\ref{homeomorphisms}.  Let us justify that the space $Z$ is well defined.

\begin{Lemma}\label{equivZ}
Suppose that $(x,\p)\sim (y,\p)$ and $(y,\p)\sim (z,\p)$. Then there exists $w\leq x,y,z$ with $\p\in D([w,w])$.  In particular, $\sim$ is an equivalence relation.
\end{Lemma}
\begin{proof}
Reflexivity and symmetry of $\sim$ are immediate.  Suppose $(x,\p)\sim (y,\p)$ and $(y,\p)\sim (z,\p)$.  Then we can find $u,v\in X$ so that $u\leq x,y$ and $v\leq y,z$ and $\p\in D([u,u])\cap D([v,v]) = D([u,u][v,v])$.
The dual of Lemma~\ref{meetsinetale} implies that $w=u[v,v]$ is the meet of $u$ and $v$ in $X$, and hence $w\leq x,y,z$.  Now $\p([w,w]) = \p([u[v,v],u[v,v]])= \p([v,v][u,u][v,v])=1$.  Thus $\p\in D([w,w])$ and $(x,\p)\sim (z,\p)$, completing the proof.
\end{proof}

Next we want to prove that $Z$ is locally compact and is equipped with open surjective maps to $\wh{E(S)}$ and $\wh{E(T)}$.

\begin{Prop}\label{etale}
The space $Z$ is locally compact.  Moreover, there are surjective \'etale maps $\sigma\colon Z\to \wh{E(S)}$ and $\tau\colon Z\to \wh{E(T)}$ defined by $\sigma(\ov{(x,\p)}) = x\p$ and $\tau (\ov{(x,\p)})=\p$.
\end{Prop}
\begin{proof}
First we show that the sets $(x,U)$  form a basis for a topology.  Suppose that $\ov{(z,\p)}\in (x,U)\cap (y,V)$.  By Lemma~\ref{equivZ}, there exists $w\in X$ with $w\leq x,y,z$ and $\p \in D([w,w])$.  Setting  $N=U\cap V\cap D([w,w])$, we have that $(w,N)$ is a neighborhood of $\ov{(z,\p)}$ contained in $(x,U)\cap (y,V)$.  Next, we verify that $\tau$ is a surjective \'etale map.   Let $\p\in \wh{E(T)}$ and suppose $e\in E(T)$ is such that $\p(e)=1$.  Choose $x\in X$ with $e=[x,x]$.  Then $\tau(\ov{(x,\p)}) = \p$ and so $\tau$ is onto.   Now we turn to the continuity of $\tau$.  Let $U$ be an open neighborhood of $\p=\tau(\ov{(x,\p)})$.  Then $(x,U\cap D([x,x]))$ is a neighborhood of $\ov{(x,\p)}$ with $\tau((x,U\cap D([x,x]))) \subseteq U\cap D([x,x])\subseteq U$.

Next we establish that, for any basic open subset $(x,U)$ of $Z$, one has $\tau\colon (x,U)\to U$ is a homeomorphism.
Define a map $\rho\colon U\to (x,U)$ by $\rho(\p) = \ov{(x,\p)}$.  This is well defined since $U\subseteq D([x,x])$.  It is immediate that $\tau|_{(x,U)}$ and $\rho$ are inverses.  Let $(y,W)$ be a basic open neighborhood of $\ov{(x,\p})$ contained in $(x,U)$.  Then since $\ov{(x,\p})\in (y,W)$ there exists $u\leq x,y$ with $\p([u,u]) =1$.  Then $N=W\cap D([u,u])$ is a neighborhood of $\p$ in $U$ and $\rho(N) = (x,N)=(y,N)\subseteq (y,W)$.  It follows that $\tau$ is \'etale.

Finally, we show that $\sigma$ is a surjective \'etale map.  For surjectivity, let $\p\in \wh{E(S)}$ and choose $e\in E(S)$ with $\p(e)=1$.  Choose $x\in X$ with $\langle x,x\rangle = e$.  Then $x^*\p([x,x]) =\p(\langle x[x,x],x[x,x]\rangle) = \p(\langle x,x\rangle) = \p(e)=1$, so $\ov{(x,x^*\p)}\in Z$.  Moreover, $\sigma( \ov{(x,x^*\p)}) = xx^*\p = \p$ as $\p\in D(\langle x,x\rangle)$. Therefore, $\sigma$ is surjective.  Next, we observe that if $(x,U)$ is a basic neighborhood of $Z$, then $\sigma|_{(x,U)} =\alpha_x\tau|_{(x,U)}$ and hence is  a homeomorphism with its image, the open subset $\alpha_x(U)$ of $\wh{E(S)}$.  It follows $\sigma$ is \'etale.

Since $Z$ is locally homeomorphic to the locally compact Hausdorff space $\wh{E(T)}$, this completes the proof.
\end{proof}

Our  ultimate goal is to show that $Z,\sigma,\tau$ give a Morita equivalence between $\mathscr G(S)$ and $\mathscr G(T)$.  First we need to decongest notation.  Let's start with $\mathscr G(T)[Z]$.  A typical element looks like $\gamma = (\ov{(x',\p)},\ov{(t,\psi)},\ov{(x,\rho)})$ with $x,x'\in X$, $\p\in D([x',x'])$, $\psi\in D(t^*t)$, $\rho\in D([x,x])$ and $\rho=\tau(\ov{(x,\rho)}) = \psi$, $\p= \tau(\ov{(x',\p)})=t\psi$.  Thus we can effectively drop $\rho$ and $\p$ from the notation and write elements in the form $\ov{(x',t,\psi,x)}$ where $\psi\in D(t^*t)\cap D([x,x])$ and $t\psi\in D([x',x'])$.  Let us determine when two quadruples $(x',t,\p,x)$ and $(y',u,\psi,y)$ represent the same element of $\mathscr G(T)[Z]$.  In order, for this to happen $\p = \psi$ and there must be $z',z\in X$ and $v\in T$ with $z'\leq x',y'$, $z\leq x,y$ and $v\leq t,u$ with $\p\in D(v^*v)\cap D([z,z])$ and $v\p\in D([z',z'])$.  In particular, if $\gamma=\ov{(x',t,\p,x)}$ and $\rho=\ov{(x',u,\p,x)}$ are elements of $\mathscr G(T)[Z]$ with $t\leq u$, then $\gamma=\rho$.

Let  $(x',U')\times (t,V)\times (x,U)$ be a basic neighborhood of $\gamma$.  Set $W= U\cap V\cap t^*U'$.  Then $(x',tW)\times (t,W)\times (x,W)$ is a smaller neighborhood of $\gamma$, which we denote by $(x',t,W,x)$.  Note that
\begin{equation}\label{wconditions}
W\subseteq D(t^*t)\cap D([x,x])\ \text{and}\ tW\subseteq D([x',x']).
\end{equation}
In terms of our decongested notation, $(x',t,W,x) = \{\ov{(x',t,\p,x)}\mid \p\in W\}$.  In summary, a basis for the topology of $\mathscr G(T)[Z]$ is given by all subsets of the form $(x',t,W,x)$ where $W$ is an open set of $\wh{E(T)}$ such that \eqref{wconditions} holds.  The groupoid structure in this notation is given by $d(\ov{(x',t,\p,x)}) = \ov{(x,\p)}$ and $r(\ov{(x',t,\p,x)}) =\ov{(x',t\p)}$.  The product is given by $\ov{(x'',t,\p,x')}\cdot \ov{(x',u,\psi,x)}= \ov{(x'',tu,\psi,x)}$ and the inverse by $\ov{(x',t,\p,x)}\inv = \ov{(x,t^*,t\p,x')}$.

Similarly, elements of $\mathscr G(S)[Z]$ can be represented in the form $\ov{(x',s,\p,x)}$ where $x,x'\in X$, $s\in S$, $\p\in D(s^*s)\cap D(\langle x,x\rangle)$ and $s\p \in D(\langle x',x'\rangle)$. Two quadruples $(x',s,\p,x)$ and $(y',s',\psi,y)$ are equivalent if and only if $\p=\psi$ and there exist $z',z\in X$ and $v\in S$ with $z'\leq x',y'$, $z\leq x,y$ and $v\leq s,s'$ with $\p\in D(v^*v)\cap D(\langle z,z\rangle)$ and $v\p\in D(\langle z',z'\rangle)$.  The correspondence takes $\gamma = (\ov{(x',\p)},\ov{(s,\psi)},\ov{(x,\rho)})$ to $\ov{(x',s,\psi,x)}$ using that $\rho = x^*\psi$ and $\p = (x')^*s\psi$.  Using Proposition~\ref{homeomorphisms} and an argument as above, one easily shows that the sets of the form $(x',s,U,x) = \{\ov{(x',s,\p,x)}\mid \p \in U\}$ where $U$ is an open subset of $\wh{E(S)}$ with $U\subseteq D(s^*s)\cap D(\langle x,x\rangle)$ and $sU\subseteq D(\langle x',x'\rangle)$ form a basis for the topology on $\mathscr G(S)[Z]$.  The groupoid structure is given by $d(\ov{(x',s,\p,x)}) = \ov{(x,x^*\p)}$ and $r(\ov{(x',s,\p,x)}) =\ov{(x',(x')^*s\p)}$.  The product is given by $\ov{(x'',s,\p,x')}\cdot \ov{(x',s',\psi,x)}= \ov{(x'',ss',\psi,x)}$ and the inverse by $\ov{(x',s,\p,x)}\inv = \ov{(x,s^*,s\p,x')}$.

\begin{Thm}\label{Moritagroupoids}
Let $S$ and $T$ be strongly Morita equivalent inverse semigroups.  Then their respective universal groupoids $\mathscr G(S)$ and $\mathscr G(T)$ are Morita equivalent.
\end{Thm}
\begin{proof}
Let $X$ be an equivalence bimodule and let $Z,\sigma,\tau$ be as in Proposition~\ref{etale}.   The maps $\sigma$ and $\tau$ are open, in fact \'etale, and surjective.  Let us show that they give rise to a Morita equivalence between  $\mathscr G(S)$ and $\mathscr G(T)$.  We shall use our shorthand notation for elements of $\mathscr G(S)[Z]$ and $\mathscr G(T)[Z]$ throughout the proof.  Let us define  $\Phi\colon \mathscr G(S)[Z]\to \mathscr G(T)[Z]$ and $\Psi\colon \mathscr G(T)[Z]\to \mathscr G(S)[Z]$ as the identity map on objects and on arrows by
\begin{align*}
\Phi(\ov{(x',s,\p,x)}) &= \ov{(x',[x',sx], x^*\p,x)}\\
\Psi(\ov{(x',t,\p,x)}) &= \ov{(x',\langle x't,x\rangle, x\p,x)}.
\end{align*}
Our goal is to show that $\Phi$ and $\Psi$ give a topological isomorphism of $\mathscr G(S)[Z]$ and $\mathscr G(T)[Z]$.  We begin by verifying that $\Psi$ is well defined.

First note that if $s=\langle x't,x\rangle$, then \[s^*s = \langle x,x't\rangle\langle x't,x\rangle = \langle x[x't,x't],x\rangle = \langle xt^*[x',x']t,x\rangle\] by Proposition~\ref{basicprops}.  Therefore \[s^*sx= \langle xt^*[x',x']t,x\rangle x= xt^*[x',x']t[x,x] = xt^*[x',x']t\] and hence we have, for $\p\in D([x,x])$,
\begin{align*}
x\p(s^*s) &= \p([s^*sx,s^*sx]) = \p([xt^*[x',x']t, xt^*[x',x']t])\\  &= \p(t^*[x',x']t[x,x]t^*[x',x']t) = t\p([x',x'])\p([x,x])=1.
\end{align*}
We conclude $x\p\in D(\langle x't,x\rangle^*\langle x't,x\rangle)$.  From $\p\in D([x,x])$, we obtain $x\p\in D(\langle x,x\rangle)$ by Proposition~\ref{homeomorphisms}.  On the other hand, using $t\p\in D([x',x'])$ and Lemma~\ref{explicitref} with $s=\langle x't,x\rangle$, we obtain
\begin{equation}\label{willneedagain}
\langle x't,x\rangle x\p = x't[x,x]\p = x't\p\in D(\langle x',x'\rangle)
\end{equation}
where the last equality uses $[x,x]\p=\p$ as $\p\in D([x,x])$.  We conclude $\ov{(x',\langle x't,x\rangle, x\p,x)}\in \mathscr G(S)[Z]$.
 Suppose $(x',t,\p,x)$ and $(y',u,\p,y)$ determine the same element of $\mathscr G(T)[Z]$. So there exist $z',z\in X$ and $v\in T$ so that: $\p\in D(v^*v)\cap D([z,z])$; $v\p\in D([z',z'])$; $z'\leq x',y'$; $z\leq x,y$; and \mbox{$v\leq t,u$}.    By Proposition~\ref{etalestruct}, $\langle z'v,z\rangle\leq \langle x't,x\rangle, \langle y'u,y\rangle$.  By Proposition~\ref{restriction}, $z\p=x\p = y\p$.  Thus $\ov{(x',\langle x't,x\rangle, x\p,x)}= \ov{(y',\langle y'u,y\rangle, y\p,y)}$, as required.  Similarly, one verifies that $\Phi$ is well defined.

Notice that both groupoids have unit space $Z$.  Let us show that $\Psi$ is a functor.   Let $\gamma =   \ov{(x',t,\p,x)}\in \mathscr G(T)[Z]$.  Then \[d(\Psi(\gamma)) = d(\ov{(x',\langle x't,x\rangle, x\p,x)}) = \ov{(x,x^*x\p)} = \ov{(x,\p)} =\Psi(d(\gamma)),\] whereas $r(\Psi(\gamma)) = r(\ov{(x',\langle x't,x\rangle, x\p,x)}) = \ov{(x',(x')^*\langle x't,x\rangle x\p)} =\ov{(x',t\p)}=\Psi(r(\gamma))$ by \eqref{willneedagain}.  Next we verify that $\Psi$ preserves products.  So suppose that $\gamma= \ov{(x'',t,\p,x')}$ and $\lambda = \ov{(x',u,\psi,x)}$ are composable, whence $\p = u\psi$.  Then we have
\begin{align*}
\Psi(\gamma)\Psi(\lambda) &=  \ov{(x'',\langle x''t,x'\rangle, x'\p,x')}\cdot \ov{(x',\langle x'u,x\rangle, x\psi,x)}\\ &= \ov{(x'',\langle x''t,x'\rangle\langle x'u,x\rangle,x\psi,x)}\end{align*}
whereas $\Psi(\gamma\lambda) = \Psi(\ov{(x'',tu, \psi,x)}) =  \ov{(x'',\langle x''tu,x\rangle, x\psi,x)}$.  Now we compute \[\langle x''t,x'\rangle\langle x'u,x\rangle = \langle x''t[x',x'u],x\rangle = \langle x''t[x',x']u,x\rangle\leq \langle x''tu,x\rangle\] by Proposition~\ref{etalestruct}, as $x''t[x',x']u\leq x''tu$.  It follows that \[\ov{(x'',\langle x''t,x'\rangle\langle x'u,x\rangle,x\psi,x)}=\ov{(x'',\langle x''tu,x\rangle, x\psi,x)}\] and hence $\Psi(\gamma)\Psi(\lambda) = \Psi(\gamma\lambda)$.  The proof that $\Phi$ is a functor is symmetric and so we leave it to the reader.

We now establish that $\Phi$ and $\Psi$ are inverses.  Let $\gamma =  \ov{(x',t,\p,x)}\in \mathscr G(T)[Z]$.  Let $s= \langle x't,x\rangle\in S$. It is straightforward to compute
\[\Phi\Psi(\gamma) = \Phi(\ov{(x',s, x\p,x)}) = \ov{(x',[x',sx],x^*x\p,x)} = \ov{(x',[x',sx],\p,x)}\] where the last equality uses Proposition~\ref{homeomorphisms}.  Now $sx = \langle x't,x\rangle x = x't[x,x]$  and so \[[x',sx] = [x',x't[x,x]] = [x',x']t[x,x]\leq t.\]  It follows that $\ov{(x',[x',sx],\p,x)}=\ov{(x',t,\p,x)}$ and so $\Phi\Psi(\gamma) = \gamma$.  An entirely symmetric argument shows that $\Psi\Phi$ is the identity and hence
$\Phi$ and $\Psi$ are inverses.
It remains to prove that they are continuous.  To do so, it suffices to show that they both send basis sets to open sets.
We just handle $\Psi$, as the case of $\Phi$ is completely analogous.

Let $(x',t,W,x)$ be a basic neighborhood of $\mathscr G(T)[Z]$; so \eqref{wconditions} holds.
We claim that $\Psi(x',t,W,x)=(x',\langle x't,x\rangle,xW,x)$ and that the latter set is open.  First we verify that the right hand side is a basic neighborhood of $\mathscr G(S)[Z]$.  The set $xW$ is open by Proposition~\ref{homeomorphisms}.  Setting $s=\langle x't,x\rangle$, we need to show that $xW\subseteq D(s^*s)\cap D(\langle x,x\rangle)$ and $s(xW)\subseteq D(\langle x',x'\rangle)$.   But this is immediate from the proof that $\Psi$ is well defined.
The inclusion $\Psi(x',t,W,x)\subseteq (x',\langle x't,x\rangle,xW,x)$ is obvious.  For the reverse inclusion, we again set $s=\langle x't,x\rangle$ and suppose that $\ov{(x',s,\p,x)}\in (x',s,xW,x)$.  Then $x^*\p\in x^*xW=W$.  Equation \eqref{wconditions} implies that $\ov{(x',t,x^*\p,x)}\in (x',t,W,x)$ and clearly $\Psi(\ov{(x',t,x^*\p,x)})=\ov{(x',s,\p,x)}$.  This completes the proof that $\mathscr G(S)$ and $\mathscr G(T)$ are Morita equivalent.
\end{proof}

In~\cite{Exel}, Exel associates to an inverse semigroup $S$ with zero, an \'etale groupoid $\mathscr G(S)_{tight}$ by taking the reduction of $\mathscr G(S)$ to the closure of the subspace of ultrafilters in $\wh{E(S)}$, i.e., taking the full subgroupoid whose unit space is the closure of the space of ultrafilters.  Recall that a homomorphism $\p\colon E(S)\to \{0,1\}$ is determined by the filter $\p\inv(1)$ and that an ultrafilter is a maximal proper filter.   The reader can easily verify that if $S$ and $T$ are strongly Morita equivalent inverse semigroups with zero, then the Morita equivalence in Theorem~\ref{Moritagroupoids} restricts to a Morita equivalence between $\mathscr G(S)_{tight}$ and $\mathscr G(T)_{tight}$.

As a corollary we obtain the strong Morita equivalence of the $C^*$-algebras of strongly Morita equivalent inverse semigroups.

\begin{Cor}\label{algebraequivalence}
Let $S$ and $T$ be strongly Morita equivalent inverse semigroups.  Then the universal and reduced $C^*$-algebras of $S$ and $T$ are strongly Morita equivalent.
\end{Cor}

In the next section we shall provide examples showing that the converse of Theorem~\ref{Moritagroupoids} is false.

Recall that $S$ is called an \emph{$F$-inverse semigroup} if $\sigma^{-1}(g)\cap eSf$ has a maximum element for all $g\in G(S)$ and $e,f\in E(S)$ where $\sigma\colon S\to G(S)$ is the maximal group image homomorphism.  An $F$-inverse semigroup is automatically $E$-unitary.  They were characterized by Lawson in terms of their representation as a $P$-semigroup~\cite{Lawson}.  Let us remind the reader of this construction.

 A \emph{McAlister triple} $(X,Y,G)$ consists of a group $G$ acting on a poset $X$ containing a downset $Y$ of $X$ so that: $Y$ has all binary meets, $G\cdot Y=X$ and $gY\cap Y\neq \emptyset$ for all $g\in G$.  The associated \emph{$P$-semigroup} $P(X,Y,G)$ consists of all pairs $(y,g)\in Y\times G$ so that $g\inv y\in Y$.  Multiplication is given by \[(y,g)(y',g') = (g(g\inv y\wedge y'),gg') = (y\wedge gy',gg')\] and inversion is given by $(y,g)^* = (g\inv y,g\inv)$.  A $P$-semigroup is $E$-unitary with semilattice of idempotents $Y$ and maximal group image $G$.  Conversely, McAlister proved that if $S$ is any $E$-unitary inverse semigroup, then $S\cong P(X,Y,G)$ for a unique McAlister triple (up to isomorphism of triples)~\cite{McAlisterPstuff,Lawson}.  Of course $G=G(S)$ and $Y=E(S)$.  The set $X$ can be obtained as the globalization~\cite{Lawsonkellendonk} (or enveloping action~\cite{abadie}) of the partial action of $G$ on $Y=E(S)$~\cite{Lawsonkellendonk}.  More concretely~\cite{topos}, one can describe $X$ as the poset of cyclic subsets (ordered by inclusion) of the \'etale right $S$-set $(G\times E(S),p)$ where $(g,e)s = (g\sigma(s),s^*es)$ and $p((g,e))=e$.  By a cyclic subset, we mean a subset of the form $(g,e)S$.   For other approaches, see~\cite{Lawson}.

 Suppose that $S\cong P(X,Y,G)$ is an $E$-unitary inverse semigroup.  Then Lawson proved that $S$ is $F$-inverse if and only if $X$ is a meet semilattice.  Moreover, in this case the semidirect product $X\rtimes G$ is an enlargement of $S$~\cite{Lawson} and hence $S$ is strongly Morita equivalent to $X\rtimes G$.   In fact, Lawson showed that $F$-inverse semigroups are precisely those semigroups admitting a semidirect product $X\rtimes G$ with $G$ a group and $X$ a semilattice as an enlargement.
Khoskam and Skandalis~\cite{Skandalis} proved as a special case of a more general result that the universal and reduced $C^*$-algebras of an $F$-inverse semigroup $S$ are strongly Morita equivalent to, respectively, the universal and reduced cross products of $G(S)$ with a certain commutative $C^*$-algebra.  We give a new proof of their result from a semigroup viewpoint that makes the commutative $C^*$-algebra in question very explicit.  We first recall a standard construction from groupoid theory.

Let $G$ be a group acting on a locally compact Hausdorff space $X$.  Then associated transformation groupoid $G\ltimes X$ has underlying topological space $G\times X$ and unit space $X$.  The groupoid structure is given by $d(g,x) = x$, $r(g,x) = gx$, $(g,x)(h,y) = (gh,y)$ and $(g,x)\inv = (g\inv,gx)$.  One easily verifies that $G\ltimes X$ is a Hausdorff \'etale groupoid.  Moreover, it is well known that $C^*(G\ltimes X)\cong C_0(X)\rtimes G$ and $C^*_r(G\ltimes X) = C_0(X)\rtimes_r G$~\cite{Renault}.  Our next proposition shows that if $E$ is a semilattice and $G$ is a group acting on $E$, then the universal groupoid $\mathscr G(E\rtimes G)$ is topologically isomorphic to $G\ltimes \wh{E}$.

\begin{Prop}\label{iscrossproduct}
Let $E$ be a semilattice and $G$ a group acting on $E$ by automorphisms.  Then $\mathscr G(E\rtimes G)\cong G\ltimes \wh{E}$ and consequently $C^*(E\rtimes G)\cong C_0(\wh E)\rtimes G\cong C^*(E)\rtimes G$ and $C^*_r(E\rtimes G)\cong C_0(\wh E)\rtimes_r G\cong  C^*(E)\rtimes_r G$.
\end{Prop}
\begin{proof}
The action of $G$ on $E$ induces an action of $G$ on $\wh E$ by $(g\p)(e) = \p(g\inv e)$.  Note that the natural partial order on $E\rtimes G$ is given by the product order, i.e., $(e,g)\leq (f,h)$ if and only if $g=h$ and $e\leq f$.   The idempotent set of $E\rtimes G$ is $E\times \{1\}\cong E$, which we identify with $E$ from now on.  A formula we shall use frequently is  $(e,g)^*(e,g) = (g\inv e,1)$.  Suppose that $\p\in D(g\inv e)\cap D(g\inv{f})$.   Then $\ov{((e,g),\p)}=\ov{((f,g),\p)}$ since $\p\in D(g\inv{(ef)})$ and $(ef,g)\leq (e,g),(f,g)$.  Thus there is a well defined surjective map \mbox{$\Phi\colon G\ltimes \wh{E}\to \mathscr G(E\rtimes G)$} given by $\Phi(g,\p) = \ov{((e,g),\p)}$ where $e\in E$ is any element so that $\p(g\inv e)=1$.  Define $\Phi$ to be the identity on objects.  Note that $\Phi$ is injective, since $\Phi(g,\p)=\Phi(h,\psi)$ implies $\ov{((e,g),\p)} = \ov{((f,h),\psi)}$ where $\p(g\inv e)=1=\psi(h\inv f)$.  But then $\p =\psi$ and $(e,g),(f,h)$ have a common lower bound, which implies $g=h$.  We proceed to show that $\Phi$ is a functor.

Clearly $\Phi$ respects domains.  Suppose $\p(g\inv e)=1$.  Then $\Phi(g,\p) = \ov{((e,g),\p)}$ and $r(g,\p) = g\p$, whereas $r \ov{((e,g),\p)} = (e,g)\p$.  Now if $f\in E$, then on the one hand $(g\p)(f) = \p(g\inv f)$; on the other hand since $(e,g)^*(f,1)(e,g) = (g\inv e,g\inv)(f,1)(e,g)= (g\inv(efe),1)$, it follows that
$(e,g)\p(f) = \p(g\inv(ef)) = \p(g\inv f)$ since $\p(g\inv e)=1$.  Thus $\Phi r(g,\p)=g\p =(e,g)\p= r\Phi(g,\p)$.  It remains to show that $\Phi$ preserves the product.  Note that if $g\p=\psi$, then $(h,\psi)(g,\p) = (hg,\p)$.  Choose $e,f\in E$ so that $\p(g\inv e)=1$ and $\psi(h\inv f)=1$.  Then by the above $(e,g)\p = g\p = \psi$ and so in $\mathscr G(E\rtimes G)$ we have $\ov{((f,h),\psi)}\cdot\ov{((e,g),\p)} = \ov{((f\cdot he,hg),\p)}=\Phi(hg,\p)$.  The last equality follows since \[\p(g\inv h\inv (f\cdot he)) = \p(g\inv e)\p(g\inv h\inv f) =\p(g\inv e)\psi(h\inv f)=1.\]  This establishes that $\Phi$ is an isomorphism of groupoids.

Now we need to check that $\Phi$ is continuous and open.  Let $\{g\}\times U$ be a basic neighborhood of $G\ltimes \wh E$ and $\p\in U$.  Choose $e\in E$ so that $\p(g\inv e)=1$.  Then $\Phi(g,\p) = \ov{((e,g),\p)}$. We claim that $((e,g),U)$  is a neighborhood of $\Phi(g,\p)$ contained in $\Phi(\{g\}\times U)$.  Indeed, if $\ov{((e,g),\psi)}\in ((e,g),U)$, then $(g,\psi)\in \{g\}\times U$.  Moreover, $\psi(g\inv e)=1$ implies that $\Phi(g,\psi) = \ov{((e,g),\psi)}$.   This shows that $\Phi$ is an open map.  Next we turn to continuity.  Let $((f,h),U)$ be a neighborhood of $\Phi(g,\p)$.  Suppose that $\p(g\inv e)=1$.  Then $\Phi(g,\p) = \ov{((e,g),\p)}\in ((f,h),U)$.  Thus $\p\in U$ and $g=h$ since $(e,g)$ and $(f,h)$ have a common lower bound.   Also $U\subseteq D(g\inv f)$.  Thus $\p(g\inv f)=1$.    The set $\{g\}\times U$ is a neighborhood of $(g,\p)$.  We claim $\Phi(\{g\}\times U)\subseteq ((f,g),U)$.  Indeed, if $\psi\in U$, then $\Phi(g,\psi) = \ov{((f,g),\psi)}\in ((f,g),U)$.   Thus $\Phi$ is a homeomorphism, completing the proof of the theorem.
\end{proof}

We now apply the above result to $F$-inverse semigroups.

\begin{Thm}
Let $S$ be an $F$-inverse semigroup with $P$-representation $S=P(X,Y,G)$; in particular, $Y=E(S)$ and $G=G(S)$.   Then $C^*(S)$ and $C^*_r(S)$ are strongly Morita equivalent to the universal cross product algebra $C^*(X)\rtimes G\cong C_0(\wh X)\rtimes G$ and the reduced cross product algebra $C^*(X)\rtimes_r G\cong C_0(\wh X)\rtimes_r G$, respectively.
\end{Thm}
\begin{proof}
The semidirect product $X\rtimes G$ is an enlargement of $S$ and hence they are strongly Morita equivalent.   By Proposition~\ref{iscrossproduct}, $C^*(X\rtimes G)\cong C^*(X)\rtimes G$ and $C^*_r(X\rtimes G)\cong C^*(X)\rtimes_r G$.   Corollary~\ref{algebraequivalence} yields the desired result.
\end{proof}

Let us present some examples.

\begin{Example}[Birget-Rhodes/Exel expansion]
Let $G$ be a group and denote by $P_{fin}(G)$ the semilattice of finite non-empty subsets of $G$ ordered by reverse inclusion (so union is the binary operation).  Then $P_{fin}(G)$ is the free semilattice on $G$.  The group $G$ acts naturally on $X=P_{fin}(G)$.  The set $Y=\{A\in X\mid 1\in A\}$ is a downset in $X$ such that $G\cdot Y=X$ and $gY\cap Y\neq \emptyset$ for all $g\in G$ (consider $\{1,g\inv\}$).  Thus we can form the $P$-semigroup $G^{\mathsf{Pr}}=P(X,Y,G)$, which is an $F$-inverse monoid known as the Birget-Rhodes expansion of $G$~\cite{BR--exp,birgetrhodesexp,InvExp}.  It was shown by Kellendonk and Lawson~\cite{Lawsonkellendonk} that $G^{\mathsf{Pr}}$ is isomorphic to the inverse semigroup associated to $G$ by Exel in his study of partial group actions~\cite{exelpartial}. In fact, the author, together with Lawson and Margolis, showed directly~\cite{InvExp} that Exel's universal property for $G^{\mathsf{Pr}}$~\cite{exelpartial} is equivalent to Szendrei's~\cite{birgetrhodesexp}.
Since $P_{fin}(G)$ is a free semilattice on $G$, it is easy to see that $\wh{P_{fin}(G)}\cong \{0,1\}^{G}\setminus \{\ov 0\}$ where $\ov 0$ is the constant function to $0$.  The action of $G$ on $\wh{P_{fin}(G)}$ is the restriction of the usual Bernoulli shift action.  So the $C^*$-algebra of $G^{\mathsf {Pr}}$ is strongly Morita equivalent to $C_0(\{0,1\}^G\setminus \{\ov 0\})\rtimes G$.
\end{Example}

\begin{Example}[free inverse monoids]  Let $A$ be a set and $FIM(A)$ the free inverse monoid generated by $A$.  Then $FIM(A)$ is an $F$-inverse monoid with maximal group image the free group $F_A$ on $A$~\cite{Lawson}.  A McAlister triple for $FIM(A)$ is $(X,Y,F_A)$ where $Y$ is the set of all finite subtrees of the Cayley graph $T$ of $F_A$ containing the vertex $1$ and $X$ is the set of all non-empty finite subtrees of $F_A$~\cite{Lawson}. Here $X$ is ordered by reverse inclusion.  If $T_1,T_2\in X$ are subtrees of $T$, then $T_1$ and $T_2$ have a meet in $X$.  More precisely, $T_1\wedge T_2$ is the subtree spanned by $T_1$ and $T_2$ (so it consists of $T_1\cup T_2$ and the ``bridge'' between them).

It is not hard to see that filters on $X$ are in bijection with non-empty subtrees of $T$.  Given a subtree $T'$ of $T$, one has the filter of all non-empty finite subtrees of $T$ contained in $T'$.  Conversely, if $\mathscr F$ is a filter on $X$, then the poset $(\mathscr F,\subseteq)$ is directed and so it is easy to see that $T'=\bigcup \mathscr F$ is a subtree of $T$ and $T_0\in \mathscr F$ if and only if $T_0\subseteq T'$.  The topology on $\wh X$ (viewed as the space of subtrees) is the subspace topology of $\{0,1\}^{V(T)\cup E(T)}$ where $V(T)$ is the vertex set of $T$ and $E(T)$ is the edge set of $T$.
The action of $F_A$ on $\wh X$ is by translation of subtrees.  So $C^*(FIM(A))$ is strongly Morita equivalent to $C_0(\wh X)\rtimes F_A$.
\end{Example}

It is much easier to show in the discrete setting that strongly Morita equivalent inverse semigroups have Morita equivalent algebras.
\begin{Thm}
Let $S$ and $T$ be strongly Morita equivalent inverse semigroups and let $K$ be a commutative ring with unit.  Then the semigroup algebras $KS$ and $KT$ are Morita equivalent.
\end{Thm}
\begin{proof}
Let  $(S,T,X,\langle\ ,\ \rangle,[\ ,\ ])$  be a Morita context.  Then $KX$ is naturally a $KS$-$KT$ bimodule.  It is also a $KT$-$KS$ bimodule via the involutions.  The map $x\to \langle x,x\rangle$ extends to a surjective $KT$-bilinear map $KX\times KX\to KS$ and similarly $[\ ,\ ]$ extends to a $KS$-bilinear map making $(KS,KT,KX,KX,\langle\ ,\ \rangle,[\ ,\ ])$ a Morita context in the ring theoretic sense.  Since $KS$ and $KT$ are rings with local units~\cite{localunits}, this establishes the Morita equivalence of $KS$ and $KT$.
\end{proof}

\section{Properties invariant under strong Morita equivalence}
This section is devoted to elaborating on the many structural properties shared by strongly Morita equivalent inverse semigroups. They have isomorphic maximal subgroups, the same maximal group image, the same set of $\D$-classes, isomorphic lattices of two-sided ideals, the same maximal group image, isomorphic cohomology groups and equivalent classifying toposes.

To expedite a proof of the above properties, we introduce two categories associated with an inverse semigroup $S$ and interpret structural properties of  $S$ in terms of these categories.  Then we shall show that the categories associated to strongly Morita equivalent inverse semigroups are naturally equivalent.  The reader is referred to~\cite{Mac-CWM} for basics on category theory.  See also~\cite{MM-Sheaves} for an introduction to topos theoretic notions.
First we define a category $S_E$ with object set $E(S)$ and arrow set $\{(f,s,e)\mid e,f\in E(S), s\in fSe\}$.  An arrow $(f,s,e)$ has domain $e$ and range $f$.  The composition is given by $(f,s,e)(e,s',d) = (f,ss',d)$.  The identity at $e$ is the arrow $(e,e,e)$.  Notice the endomorphism monoid at $e$ is isomorphic to $eSe$ via the map $(e,s,e)\mapsto s$.  Hence the automorphism group at $e$ is isomorphic to the maximal subgroup $G_e$ at $e$.    The category $S_E$ is known as the \emph{idempotent splitting} (or Cauchy completion or Karoubi envelope) of $S$.     Also one can show that the isomorphisms of $S_E$ are precisely the arrows of the form $(ss^*,s,s^*s)$ (with inverse $(s^*s,s^*,ss^*)$)~\cite{funktopos}.  In fact, the groupoid of isomorphisms of $S_E$ is precisely the underlying groupoid of $S$~\cite{Lawson,funktopos}.  Recall that idempotents $e,f\in E(S)$ are $\D$-equivalent, written $e\D f$, if there exists $s\in S$ such that $ss^*=e$, $s^*s=f$~\cite{Lawson}.  Notice that two idempotents are $\D$-equivalent if and only if they are isomorphic as objects of $S_E$.

An important subcategory of $S_E$, first introduced by Loganathan in his study of inverse semigroup cohomology~\cite{loganathan}, is $L(S)$, which consists of all the objects of $S_E$ and of those arrows of the form $(f,s,e)$ such that $s^*s=e$.  In particular, it contains all the isomorphisms of $S_E$.   It is easy to verify~\cite{funktopos} that $L(S)$ consists precisely of the split monomorphisms in $S_E$; recall that a monomorphism $f\colon c\to d$ in a category is said to be \emph{split} if there exists $g\colon d\to c$ so that $gf=1$.  In particular, if $F\colon S_E\to T_E$ is an equivalence of categories, then $F$ restricts to an equivalence of categories $F\colon L(S)\to L(T)$.

The category $L(S)$ is a category of monics with pullbacks.  In the case that $S$ is semilattice, $L(S)$ is the usual category associated to the poset $S$.   It is known that the classifying topos $\mathscr BS$ of all \'etale right $S$-sets is equivalent to the category $\pv{Set}^{L(S)^{op}}$ of presheaves on $L(S)$ (contravariant functors from $L(S)$ to the category of sets)~\cite{funktopos,topos,etendues}.  In particular, if $L(S)$ is equivalent to $L(T)$, then $\mathscr BS$ is equivalent to $\mathscr BT$~\cite{funktopos}.  By Loganathan's results~\cite{loganathan}, the cohomology of $S$ in the sense of Lausch~\cite{lausch} is the cohomology of $\mathscr BS$ in the sense of Grothendieck and so if $L(S)$ is equivalent to $L(T)$, then $S$ and $T$ have the same cohomology groups.

Let's consider an example.  Let $S$ be the polycyclic inverse monoid (or Cuntz semigroup) on a set $X$ of cardinality at least $2$.  It is presented as an inverse monoid with zero by $\langle X\mid x^*y = \delta_{x,y},\ \text{for all}\ x,y\in X\rangle$.  The non-zero elements of $S$ can be represented uniquely in the form $uv^*$ with $u,v$ in the free monoid on $X$.  The non-zero idempotents are of the form $uu^*$.  Notice that endomorphism monoid of $1$ in $L(S)$ is the free monoid on $X$ since the only arrows from $1$ to $1$ are of the form $(1,u,1)$ with $u$ in the free monoid on $X$.  Also each non-zero idempotent is $\D$-equivalent to $1$, so all non-zero idempotents are isomorphic to $1$.  Indeed, $(uu^*,u,1)$ is an isomorphism from $1$ to $uu^*$.  So $L(S)$ is equivalent to the free monoid on $X$ (viewed as a one object category) with an adjoined initial object (coming from the zero).

More generally, if $\Gamma$ is a graph and $S$ is the associated graph inverse semigroup~\cite{Paterson,graphinverse}, then $L(S)$ is equivalent to the free category on $\Gamma$ with an adjoined initial object (see~\cite{Mac-CWM} for the notion of a free category on a graph).  Indeed, the full subcategory whose objects are the empty paths is isomorphic to the free category on $\Gamma$ and contains an element from each isomorphism class of objects except $0$.  It should be clear from these examples, that one should really have a separate theory for inverse semigroups with zero.  See~\cite{Lawsonzero,Lawson0eunit} for the appropriate theory in this context.

One can define a preorder on $E(S)$ by $e\leq_{\J} f$ if $SeS\subseteq SfS$.  We write $e\J f$ if $SeS=SfS$.  It is known $e\leq_{\J} f$ if and only if there exists $x\in E(S)$ so that $e\D x$ and $x\leq f$~\cite{Lawson}.  In particular, $e\D f$ implies $e\J f$.  The lattice of downsets in $(E(S)/{\J},\leq_{\J})$ is isomorphic to the lattice of two-sided ideals in $S$.  If $C$ is any category, one can define a preorder on the objects of $C$ by $e\preceq f$ if there is an arrow from $e$ to $f$.   One can easily verify that the associated poset is an invariant of $C$ up to natural equivalence, i.e., naturally equivalent categories have isomorphic posets.
Notice that in $L(S)$ there is an arrow from $e$ to $f$ if and only if there is an element $s\in S$ with $fs=s$ and $s^*s=e$.  If such an $s$ exists then $e\D ss^*\leq f$, so $e\preceq f$ implies $e\leq_{\J} f$.  Conversely, if $e\leq _{\J} f$, then there is an element $s\in S$ with $s^*s=e$ and $ss^*\leq f$, whence $fs=s$.  Therefore, $(f,s,e)$ is an arrow from $e$ to $f$ in $L(S)$, establishing that $e\preceq f$.  Thus the poset associated to $L(S)$ is $E(S)/{\J}$.  Therefore, if $L(S)$ and $L(T)$ are equivalent, then $S$ and $T$ have isomorphic $\J$-orders and hence isomorphic lattices of two-sided ideals.

Note that $S$ has a zero element if and only if $L(S)$ has an initial object.  One can show~\cite{topos} that the maximal group image $G(S)$ is determined by $\mathscr BS$ up to equivalence and hence by $L(S)$ up to equivalence.

\begin{Thm}\label{Lcatsareequiv}
Let $S$ and $T$ be strongly Morita equivalent inverse semigroups.  Then there is an equivalence of categories $F\colon S_E\to T_E$.
\end{Thm}
\begin{proof}
Let $X$ be an equivalence bimodule.  Fix, for each idempotent $e\in E(S)$, an element $x_e\in X$ with $\langle x_e,x_e\rangle =e$.  Define $F\colon S_E\to T_E$ on objects by $F(e) = [x_e,x_e]$ and on arrows by $F(f,s,e) = (F(f),[x_f,sx_e],F(e))$.  We will show that $F$ is an equivalence of categories by showing that it is a full and faithful functor and moreover each object of $T_E$ is isomorphic to an object $F(e)$ with $e\in E(S)$.

First we verify that $F$ is a functor.   The map $F$ is well-defined since
\[[x_f,x_f][x_f,sx_e][x_e,x_e] = [x_f,\langle x_f,x_f\rangle sx_e][x_e,x_e] = [x_f, fsx_e]=[x_f,sx_e]\]
where we used that $fs=s$.
Since $[x_e,ex_e] = [x_e,\langle x_e,x_e\rangle x_e]=[x_e,x_e]$, we see that $F$ takes identities to identities.  Verifying functoriality amounts to showing that  if $s\in gSf$ and $s'\in fSe$, then $[x_g,sx_f][x_f,s'x_e] = [x_g,ss'x_e]$.  We compute \[[x_g,sx_f][x_f,s'x_e] = [x_g,\langle sx_f,x_f\rangle s'x_e] =[x_g,s\langle x_f,x_f\rangle s'x_e] = [x_g,ss'x_e]\] since $s=sf=s\langle x_f,x_f\rangle$.

Now we must check $F$ is full and faithful, beginning with faithfulness.  Suppose $F(f,s,e) = F(f,s',e)$ with $s,s'\in fSe$.  Then $[x_f,sx_e]=[x_f,s'x_e]$.  Therefore, \[s=fse=\langle x_f,x_f\rangle s\langle x_e,x_e\rangle = \langle x_f[x_f,sx_e],x_e\rangle = \langle x_f[x_f,s'x_e],x_e\rangle =s'\] as required.   Next, suppose that $t\in F(f)TF(e)$ with $e,f\in E(S)$.  We need that $t = [x_f,sx_e]$ some $s\in fSe$.  Let $s = \langle x_ft,x_e\rangle$.  First we verify that $fse=s$.  Indeed, using Proposition~\ref{basicprops}
\begin{align*}
fse &= \langle x_f,x_f\rangle \langle  x_ft,x_e\rangle \langle x_e,x_e\rangle = \langle x_ft[x_e,x_e],x_e\rangle \\ &= \langle x_ft,x_e[x_e,x_e]\rangle = \langle x_ft,x_e\rangle=s
\end{align*}
Moreover, since $[x_f,x_f]t[x_e,x_e] = F(f)tF(e)=t$ we have \[[x_f,sx_e] = [x_f,\langle x_ft,x_e\rangle x_e] = [x_f,x_ft[x_e,x_e]] = [x_f, x_f]t[x_e,x_e]=t\] as was required.

To complete the proof that $F\colon S_E\to T_E$ is an equivalence, we must show that if $e\in E(T)$, then there exists $f\in E(S)$ so that $F(f)\D e$ (since isomorphism in $T_E$ is the $\D$-relation).  Choose $x\in X$ with $[x,x]=e$.  Let $f=\langle x,x\rangle\in E(S)$.  Consider $t=[x_f,x]$.  Then \[t^*t = [x,x_f][x_f,x] = [x,\langle x_f,x_f\rangle x] =[x,\langle x,x\rangle x] = [x,x]=e\] and \[tt^*=[x_f,x][x,x_f] = [x_f,\langle x,x\rangle x_f] =  [x_f,\langle x_f,x_f\rangle x_f] = [x_f,x_f] = F(f)\] completing the proof that $F$ is an equivalence.
\end{proof}

The following structural results are immediate from Theorem~\ref{Lcatsareequiv} and the discussion before the theorem.

\begin{Cor}\label{structuralresults}
Let $S$ and $T$ be strongly Morita equivalent inverse semigroups.  Then the following hold:
\begin{enumerate}
\item The categories $S_E$ and $T_E$ are equivalent;
\item The categories $L(S)$ and $L(T)$ are equivalent;
\item For each $e\in E(S)$, there is an idempotent $f\in E(T)$ with $eSe\cong fSf$ and conversely;
\item The underlying groupoids of $S$ and $T$ are naturally equivalent;
\item There is a bijection $F\colon E(S)/{\D}\to E(T)/{\D}$ such that if $D$ is a $\D$-class of $S$ with maximal subgroup $G$, then $F(D)$ is a $\D$-class of $T$ with maximal subgroup isomorphic to $G$;
\item The posets $E(S)/{\J}$ and $E(T)/{\J}$ are isomorphic;
\item $S$ and $T$ have isomorphic lattices of two-sided ideals;
\item The classifying toposes $\mathscr BS$ and $\mathscr BT$ are equivalent;
\item $S$ and $T$ have the same cohomology groups;
\item $S$ has a zero if and only if $T$ has a zero.
\end{enumerate}
\end{Cor}

If $E$ is a semilattice, then $(E,\leq) \cong (E/{\J},\leq_{\J})$ and so we have the following corollary.

\begin{Cor}
Let $E$ and $F$ be strongly Morita equivalent semilattices.  Then $E$ is isomorphic to $F$.
\end{Cor}

\begin{Example}[Morita equivalence of groupoids vs.\ strong Morita equivalence] Notice that if $E$ is a semilattice, then $\mathscr G(E)=\wh E$ where each element of $\wh E$ is a unit.  So if $E$ and $F$ are semilattices with $\wh E\cong \wh F$, then $\mathscr G(E)\cong \mathscr G(F)$.  Now observe that the space of homomorphisms $E\to \{0,1\}$ is precisely the space of (multiplicative) homomorphisms $E\to \mathbb C$, or equivalently the space of $\mathbb C$-algebra homomorphisms $\mathbb CE\to \mathbb C$ where $\mathbb CE$ is the (discrete) semigroup algebra of $E$ (with the topology of pointwise convergence).  Now if $E$ is a finite semilattice, then Solomon proved that $\mathbb CE\cong \mathbb C^E$~\cite{Burnsidealgebra} and hence $\wh E$ depends only on the cardinality of $E$.  Thus if $E$ and $F$ are two finite semilattices of the same order, then $\wh E\cong \wh F$ although if $E\ncong F$, then $E$ and $F$ are not strongly Morita equivalent.   This shows that inverse semigroups may have Morita equivalent (or even isomorphic) universal groupoids without being strongly Morita equivalent.
\end{Example}

We are now in a position to establish converses to Proposition~\ref{enlargement} and Theorem~\ref{Lcatsareequiv} for monoids.

\begin{Cor}\label{monoidcase}
Suppose $S$ and $T$ are  inverse semigroups such that $T$ is a monoid.   Then the following are equivalent:
\begin{enumerate}
\item There exists an idempotent \mbox{$e\in E(S)$} such that $S=SeS$ and $T\cong eSe$;
\item $S$ and $T$ are strongly Morita equivalent;
\item The categories $S_E$ and $T_E$ are equivalent.
\end{enumerate}
Consequently, if $S$ and $T$ are strongly Morita equivalent with $T$ a monoid, then $S$ is an enlargement of  $T$.  In particular, if $S$ and $T$ are inverse monoids, then $S$ is strongly Morita equivalent to $T$ if and only if there exist $e\in E(S)$ and $f\in E(T)$ so that $S=SeS$, $eSe\cong T$ and $T=TfT$ and $fTf\cong S$.
\end{Cor}
\begin{proof}
The first condition implies $S$ is an enlargement of $T$ and hence they are strongly Morita equivalent.  The implication (2) implies (3) is the content of Theorem~\ref{Lcatsareequiv}.  For the final statement, let  $F\colon T_E\to S_E$ be an equivalence of categories.    Set $e=F(1)$.  Then because $F$ is an equivalence, it induces an isomorphism on endomorphism monoids.  Thus $eSe=F(1)SF(1)\cong 1T1=T$.   Because $F$ restricts to an equivalence of $L(T)$ and $L(S)$ and $1$ belongs to the unique maximum class in the order $\preceq$ associated to $L(T)$, it follows that the class of $e$ is the unique maximum element in the ordering on $L(S)$, i.e., the $\J$-class of $e$ is the maximum element in $E(S)/{\J}$.  But this is equivalent to saying that $SeS=S$.
\end{proof}

It would be interesting to know whether the converse of Theorem~\ref{Lcatsareequiv} holds in general.

Next we give a fairly general condition under which strongly Morita equivalent inverse monoids must be isomorphic.  Recall that the bicyclic monoid is given by the presentation $\langle x\mid x^*x=1\rangle$.   An inverse monoid $S$ contains as a submonoid the bicyclic monoid if and only if it admits an element that is right invertible but not a unit, or equivalently if it admits an element that is left invertible but not a unit~\cite{Lawson,CP}.

\begin{Cor}
Let $S$ and $T$ be strongly Morita equivalent inverse monoids and suppose $S$ does not contain a copy of the bicyclic monoid.  Then $S\cong T$.
\end{Cor}
\begin{proof}
We know that there is an idempotent $e$ so that $SeS=S$ and $eSe\cong T$.  Since $S$ does not contain a copy of the bicyclic monoid, its non-units form an ideal.  Hence the equality $SeS=S$ implies that $e=1$.  Thus $S=eSe\cong T$.
\end{proof}

Every proper image of the bicyclic monoid is a group~\cite{CP}, so no residually finite inverse semigroup can contain a copy of the bicyclic monoid. (Actually it is known that the bicyclic monoid cannot embed in any compact semigroup since compact semigroups are stable and the bicyclic monoid cannot embed in any stable semigroup~\cite{qtheor}.)  Also no semigroup with central idempotents contains a copy of the bicyclic monoid since its idempotents are not central.  Hence we have the following corollary.

\begin{Cor}
Suppose that $S$ and $T$ are strongly Morita equivalent monoids such that $S$ is either:
\begin{enumerate}
\item a group;
\item commutative;
\item has central idempotents;
\item is residually finite.
\end{enumerate}
Then $S$ and $T$ are isomorphic.
\end{Cor}

Next, we show that strongly Morita equivalent inverse monoids have isomorphic centers; the corresponding result for semigroups is false, as we shall see shortly.  The conceptual way to see this is to show that if $S$ is a monoid, then $Z(S)$ is isomorphic to the endomorphism monoid of the identity functor on $S_E$ in the category of endofunctors of $S_E$ with arrows natural transformations.  We opt for a direct proof using the enlargement criterion.

\begin{Prop}\label{center}
Suppose $S$ and $T$ are strongly Morita equivalent inverse monoids.  Then $Z(S)\cong Z(T)$.
\end{Prop}
\begin{proof}
Without loss of generality, we may assume that there is an idempotent $e\in S$ so that $S=SeS$ and $T=eSe$.  Define $\p\colon Z(S)\to Z(T)$ by $\p(s)=es=ese=se$.  Trivially, $\p(s)\in Z(T)$ since $est=ets=tes$.  Also, if $s,s'\in Z(S)$, then $\p(s)\p(s') = ess'e = \p(ss')$.  Next we construct an inverse for $\p$.  Write $1=s_1es_2$ with $s_1,s_2\in S$.  Then define $\psi\colon Z(T)\to Z(S)$ by $\psi(t) = s_1ts_2$.  To check that $\psi(t)$ is central, let $s\in S$ and write $s=s'es''$.  Then $\p(t)s= s_1tes_2s'es'' = s_1es_2s'ets'' = s'tes''s_1es_2= s'es''s_1ets_2=s\p(t)$.  Next, observe that $\p\psi(t) = es_1ts_2 = es_1ets_2 = tes_1es_2 = t$ for $t\in Z(T)$ and $\psi\p(s) = s_1ess_2=ss_1es_2=s$ for $s\in Z(S)$.  Thus $\p$ is the procured after isomorphism.
\end{proof}

Lawson proved that if $S$ is an enlargement of $T$, then the congruence lattices $\mathrm{Cong}(S)$ and $\mathrm{Cong}(T)$ are isomorphic~\cite{Lawsonenlarge}, to wit he showed that every congruence on $T$ extends uniquely to $S$.  We are not able to prove or disprove in general that strongly Morita equivalent inverse semigroups have isomorphic congruence lattices.  However, if one imposes the additional condition that $E(S)$ and $E(T)$ are directed posets (e.g., if $S$ and $T$ are monoids), then we can prove it.   Recall that a congruence on a category is an equivalence relation on coterminal arrows that respects the multiplication; see~\cite{Mac-CWM}.

\begin{Prop}
Let $S$ be an inverse semigroup such that $E(S)$ is directed, i.e., for all $e,e'\in E(S)$, there exists $f\in E(S)$ so that $e,e'\leq f$.  Then the congruence lattices of $S$ and $S_E$ are isomorphic.
\end{Prop}
\begin{proof}
Define $\p\colon \mathrm{Cong}(S)\to \mathrm{Cong}(S_E)$ by setting $(e,s,f)\mathrel{\p(R)} (e,t,f)$ if and only if  $s\mathrel R t$.  Evidently, $\p(R)$ is a congruence and $R\subseteq R'$ implies $\p(R)\subseteq \p(R')$.  Suppose, conversely, that $\p(R)\subseteq \p(R')$ for congruences $R,R'$ on $S$ and assume $s\mathrel R t$.  Using that $E(S)$ is directed, we can find idempotents $e,f\in E(S)$ so that $ss^*,tt^*\leq e$ and $s^*s,t^*t\leq f$.  Then $(e,s,f)\mathrel{\p(R)} (e,t,f)$ and so $(e,s,f)\mathrel {\p(R')} (e,t,f)$ by assumption.  We conclude $s\mathrel{R'} t$.  Therefore, the map $\p$ is an order embedding.  To finish the proof, it remains to establish that $\p$ is onto.
So let $R$ be a congruence on $S_E$.
\begin{Claim}
Suppose that $s,t\in eSf\cap e'Sf'$ with $e,f,e',f'\in E(S)$.  Then we have $(e,s,f)\mathrel R (e,t,f)$ if and only if $(e',s,f')\mathrel R (e',t,f')$.
\end{Claim}
\begin{proof}[Proof of claim]
Suppose that $(e,s,f)\mathrel R (e,t,f)$.  Then $s,t\in e'eSff'$ and so
\begin{align*}
(e',s,f') &= (e',e'e,e)(e,s,f)(f,ff',f')\mathrel R (e',e'e,e)(e,t,f)(f,ff',f')\\ &=(e',t,f').
\end{align*}
The reverse implication is proved in the same fashion.
\end{proof}

It follows from the claim that we can define a relation $R'$ on $S$ by $s\mathrel {R'} t$ if there exist $e,f\in E(S)$ with $s,t\in eSf$ and $(e,s,f)\mathrel R (e,t,f)$.  Let us show that $R'$ is a congruence, beginning by verifying it is an equivalence relation.  Since $(ss^*,s,s^*s)\mathrel R (ss^*,s,s^*s)$, we conclude $s\mathrel{R'} s$.  It is clear that $R$ is symmetric.  Transitivity is the sticky bit.  Suppose that $s\mathrel R t$ and $t\mathrel R u$.  Then there are idempotents $e,e',f,f'$ so that $s,t\in eSf$, $t,u\in e'Sf'$,  $(e,s,f)\mathrel R (e,t,f)$ and $(e',t,f')\mathrel R (e',u,f')$.  Choose idempotents $e'',f''$ with $e,e'\leq e''$ and $f,f'\leq f''$.  Then by the claim, \[(e'',s,f'')\mathrel R (e'',t,f'')\mathrel R (e'',u,f'')\] and so $s\mathrel{R'} u$, as required.  Finally, suppose that $u,v\in S$ and $s\mathrel{R'} t$.  We show that $usv\mathrel{R'} utv$.  Choose idempotents $e,f$ so that $s,t\in eSf$ and $(e,s,f)\mathrel R (e,t,f)$.  Then
\begin{align*}
(uu^*,usv,v^*v)&=(uu^*,ue,e,)(e,s,f)(f,fv,v^*v)\\ &\mathrel R (uu^*,ue,e,)(e,t,f)(f,fv,v^*v)= (uu^*,utv,v^*v)
\end{align*}
and so $usv\mathrel {R'} utv$.  This completes the verification that $R'$ is a congruence.

Next we show that $R=\p(R')$.  Clearly, if $(e,s,f)\mathrel R (e,t,f)$, then $s\mathrel {R'} t$ and hence $(e,s,f)\mathrel{\p(R')} (e,t,f)$.  Suppose conversely \mbox{$(e,s,f)\mathrel{\p(R')} (e,t,f)$}.  Then $s\mathrel {R'} t$ and so there exist idempotents $e',f'$ with $s,t\in e'Sf'$ and $(e',s,f')\mathrel R (e',t,f')$.  But then the claim yields $(e,s,f)\mathrel R (e,t,f)$, as required.  This completes the proof $\p$ is onto, establishing the proposition.
\end{proof}

\begin{Cor}
Suppose that $S$ and $T$ are strongly Morita equivalent inverse semigroups with $E(S)$ and $E(T)$ directed.  Then $\mathrm{Cong}(S)\cong \mathrm{Cong(T)}$.
\end{Cor}

We do not know whether $S$ strongly Morita equivalent to $T$ with $E(T)$ directed implies $S$ is an enlargement of $T$, although it seems quite plausible.

It follows from Corollary~\ref{monoidcase} that if $S$ and $T$ are strongly Morita equivalent inverse monoids, then $S$ is $E$-unitary if and only if $T$ is $E$-unitary.  More generally, if $S$ and $T$ are strongly Morita equivalent, $S$ is $E$-unitary and $E(T)$ is directed, then $T$ is also $E$-unitary as a consequence of  Corollary~\ref{structuralresults}.  However, notice that the non-$E$-unitary five element Brandt semigroup
\[B_5=\left\{\begin{pmatrix} 1 &0\\0&0\end{pmatrix},\begin{pmatrix} 0 &1\\0&0\end{pmatrix},\begin{pmatrix} 0 &0\\1&0\end{pmatrix},\begin{pmatrix} 0 &0\\0&1\end{pmatrix},\begin{pmatrix} 0 &0\\0&0\end{pmatrix},\right\}\] is an enlargement of the two-element semilattice $\{0,1\}$, which is $E$-unitary, so being $E$-unitary is not in general an invariant of strong Morita equivalence.  Also notice that $Z(B_5)$ consists of only the zero matrix and hence is not isomorphic to the center of the two-element semilattice.  So Proposition~\ref{center} fails if we remove the assumption about monoids. It turns out that being $F$-inverse is an invariant of strong Morita equivalence.

\begin{Prop}
Suppose that $S$ and $T$ are strongly Morita equivalent inverse semigroups.  Then $S$ is $F$-inverse if and only if $T$ is $F$-inverse.
\end{Prop}
\begin{proof}
Let $S$ be any inverse semigroup.  Then $S_E$ is an inverse category (that is, a category so that for all arrows $s$, there exists a unique arrow $s^*$ so that $ss^*s=s$, $s^*ss^*=s^*$); namely $(f,s,e)^* = (e,s^*,f)$.  Inverse categories have a groupoid of germs (or maximal groupoid image) obtained by identifying coterminal arrows with a common lower bound.  Let us say that an inverse category $\mathscr C$ is $F$-inverse if each $\sigma$-class has a maximum element, where $\sigma\colon \mathscr C\to \pi_1(\mathscr C)$ is the maximal groupoid image.  Clearly, being $F$-inverse is preserved by natural equivalence of categories.  Now it follows directly from the definition that $S$ is $F$-inverse if and only if $S_E$ is $F$-inverse.  We conclude that being $F$-inverse is an invariant of strong Morita equivalence by Theorem~\ref{Lcatsareequiv}.
\end{proof}

Our final result shows that the maximal group image is an invariant of a strong Morita equivalence class.  Since $G(S)$ can be constructed from $L(S)$~\cite{topos}, this follows from Theorem~\ref{Lcatsareequiv}.  We provide here a direct proof for the convenience of the reader.

\begin{Prop}\label{actionongerms}
Let $X$ be an equivalence bimodule for $S$ and $T$.  Let $Y=\varinjlim X$ be the set of germs for $X$.  Then $Y$ is a free transitive left $G(S)$-set and a free transitive right $G(T)$-set.
\end{Prop}
\begin{proof}
We prove the result for $S$, as the case of $T$ is dual.  To establish transitivity, suppose that $[x],[y]\in Y$.  Let $s=\langle x,y\rangle$.  Then $sy = \langle x,y\rangle y = x[y,y]\leq x$.  Thus $\sigma(s)[y]=[sy]=[x]$.  Next suppose that $\sigma(s)[x]=[x]$.  Then $[sx]=[x]$, so there exists $y\leq sx,x$. Set $e=\langle x,x\rangle$.  Then $se=s\langle x,x\rangle = \langle sx,x\rangle \geq \langle y,y\rangle$ and so $\langle y,y\rangle\leq s$, whence $\sigma(s)=1$.  Thus $G(S)$ acts freely and transitively on $Y$.
\end{proof}

The following lemma is well known and straightforward to prove~\cite{CP,qtheor}.

\begin{Lemma}\label{schutlemma}
Let $G$ and $H$ be groups and $X$ a set such that $G$ acts freely and transitively on the left of $X$, $H$ acts freely and transitively on the right of $X$ and the actions commute.  Fix $x_0\in X$. Then the map $\psi\colon G\to H$ taking $g\in G$ to the unique element $h\in H$ with $gx_0=x_0h$ is an anti-isomorphism.  Consequently, $G\cong H$.
\end{Lemma}

Proposition~\ref{actionongerms} and Lemma~\ref{schutlemma} admit the following corollary.

\begin{Cor}
If $S$ and $T$ are strongly Morita equivalent inverse semigroups, they have isomorphic maximal group images.
\end{Cor}

\section{Conclusion}
We have defined a notion of strong Morita equivalence for inverse semigroups and shown that it implies Morita equivalence of the corresponding universal groupoids and $C^*$-algebras.  Strongly Morita equivalent inverse semigroups were shown to have equivalent classifying toposes and to enjoy many common structural properties.  Strong Morita equivalence generalizes Lawson's notion of enlargement.  As a consequence, we recover a result of Khoshkam and Skandalis on the $C^*$-algebra of an $F$-inverse semigroup in a more explicit form.  We do not have a pair of strongly Morita equivalent inverse semigroups that does not belong to the equivalence relation generated by all pairs $(S,T)$ with $S$ an enlargement of $T$, but we suspect they exist.

\bibliographystyle{abbrv}
\bibliography{standard2}
\end{document}